\documentclass[12pt]{amsart}
\usepackage{amsmath,amsfonts,amsthm,amscd}
\newtheorem{theorem}{Theorem}[section]

\newtheorem{proposition}[theorem]{Proposition}
\newtheorem{lemma}[theorem]{Lemma}
\newtheorem{definition}[theorem]{Definition}
\newtheorem{corollary}[theorem]{Corollary}
\theoremstyle{remark}
\newtheorem*{remark}{Remark}
\numberwithin{equation}{section}
\begin{document}
\bibliographystyle{amsalpha}
\title{Bach-flat asymptotically locally Euclidean metrics} 
\author{Gang Tian}
\address{Gang Tian\\
 Department of Mathematics \\ MIT \\ Cambridge, MA 02139}
\curraddr{Department of Mathematics, Princeton University, 
Princeton, NJ 08544}
\email{tian@math.mit.edu}
\thanks{The research of the first author was 
partially supported by NSF Grant DMS-0302744.}
\author{Jeff Viaclovsky}
\address{Jeff Viaclovsky, Department of Mathematics, MIT, Cambridge, MA 02139}
\email{jeffv@math.mit.edu}
\thanks{The research of the second author was partially 
supported by NSF Grant DMS-0202477.}
\date{October 13, 2003}
\begin{abstract}
We obtain a volume growth and curvature decay result for various classes of 
complete, noncompact Riemannian metrics in dimension 4; in particular our
method applies to anti-self-dual or K\"ahler metrics with 
zero scalar curvature, and metrics with harmonic curvature. 
Similar results were obtained for Einstein metrics in \cite{Anderson}, 
\cite{BKN}, \cite{Tian}, but our analysis differs from the Einstein 
case in that (1) we consider more generally a fourth order system 
in the metric, and (2) we do not assume any pointwise Ricci 
curvature bound.
\end{abstract}
\maketitle
\section{Introduction}
 In dimension $4$, the Euler-Lagrange equations of the 
functional 
\begin{align}
\label{functional}
\mathcal{W} : g \mapsto \int_X |W_g|^2 dV_g, 
\end{align}
where $W_g$ is the Weyl curvature tensor, are given by
\begin{align}
\label{bacheq}
B_{ij} \equiv \nabla^k \nabla^l W_{ikjl} +
\frac{1}{2}R^{kl}W_{ikjl} = 0,
\end{align}
where $W_{ijkl}$ and $R^{kl}$ are the components in 
of the Weyl and Ricci tensors, respectively \cite{Besse}, 
\cite{Derdzinski}. The tensor $B_{ij}$ is 
referred to as the Bach tensor \cite{Bach}.
Note that, in particular, metrics which are locally 
conformal to an Einstein metric are Bach-flat, 
and half conformally flat metrics are Bach-flat
\cite{Besse}. 
Half conformally flat metrics are also known 
as self-dual or anti-self-dual if $W^- = 0$
or $W^+ = 0$, respectively. 

A smooth Riemannian manifold $(X,g)$ 
is called an asymptotically locally 
Euclidean (ALE) end of order $\tau$ 
if there exists a finite subgroup $\Gamma \subset SO(4)$ 
acting freely on $\mathbf{R}^4 \setminus B(0,R)$ and a 
$C^{\infty}$ diffeomorphism 
$\Psi : X \rightarrow ( \mathbf{R}^4 \setminus B(0,R)) / \Gamma$ 
such that under this identification, 
\begin{align}
g_{ij} &= \delta_{ij} + O( r^{-\tau}),\\
\ \partial^{|k|} g_{ij} &= O(r^{-\tau - k }),
\end{align}
for any partial derivative of order $k$ as
$r \rightarrow \infty$. We say an end 
is ALE of order $0$ if we can find a coordinate system 
as above with $g_{ij} = \delta_{ij} + o(1)$, 
and $\partial^{|k|} g_{ij} = o(r^{- k })$
as $r \rightarrow \infty$. 
A complete, noncompact manifold $(X,g)$ is called ALE if $X$ can be written 
as the disjoint union of a compact set and finitely many ALE ends. 

We note that ALE spaces have 
been studied often
in the literature, for example, 
they arise naturally in the positive mass theorem and
the Yamabe Problem
\cite{Bartnik}, \cite{LeeParker}, \cite{Schoen4},\cite{SchoenYau3};
in orbifold compactness of Einstein metrics and
ALE Ricci-flat metrics
\cite{Anderson}, \cite{BKN}, \cite{CheegerTian}, 
\cite{Kronheimer}, \cite{Tian}, \cite{Nakajima}; 
and in gluing theorems for scalar-flat K\"ahler metrics 
\cite{CalderbankSinger}, \cite{KovalevSinger}.
There has been a considerable amount of research on the 
existence of anti-self-dual metrics on compact manifolds, 
we mention \cite{Poon}, \cite{LeBrun2}, \cite{Floer},
\cite{DonaldsonFriedman}, \cite{Taubes}, also see 
\cite{LeBrun5} for a nice survey and further references. 
We note that on any such manifold, by using the Green's function for 
the conformal Laplacian as a conformal factor, we obtain 
a complete non-compact scalar-flat anti-self-dual 
ALE metric. By taking a finite sum of Green's functions 
based at different points as a conformal factor, we obtain examples with 
several ends. 

We recall that for a complete, noncompact manifold, 
the Sobolev constant $C_S$ is 
defined as the best constant $C_S$ so that 
for all $f \in C^{0,1}_c(X)$,  we have
\begin{align}
\label{siq}
\Vert f \Vert_{L^{\frac{2n}{n-2}}} \leq C_S \Vert \nabla f \Vert_{L^2}. 
\end{align}
\begin{theorem} 
\label{decay}
Let $(X,g)$ be a complete, noncompact 4-dimensional Riemannian 
manifold which is Bach-flat and has zero scalar curvature.
Assume that
\begin{align}
\int_X |Rm_g|^2 dV_g < \infty, \mbox{ and } C_S < \infty, 
\end{align}
where $Rm_g$ denotes the Riemannian curvature tensor. 
Fix a base point $p \in X$, then 
\begin{align}
\label{imp2}
\underset{S(r)}{sup}|\nabla^k Rm_g| = o(r^{-2-k}),
\end{align}
as $r \rightarrow \infty$, where $S(r)$ denotes sphere of radius $r$
centered at $p$.

 Assume furthermore that 
\begin{align}
\label{a2}
b_1(X) < \infty,
\end{align}
where $b_1(X)$ denotes the first Betti number. 
Then there exists a constant $C$ (depending on g) such that 
\begin{align}
\label{imp}
Vol(B(p,r)) \leq C r^4,
\end{align}
$(X,g)$ has finitely many ends, and 
each end is ALE of order $0$. 
\end{theorem}
\begin{remark} From \cite[Theorem 1]{Carron}, there is a 
bound on the number of ends depending only upon the Sobolev constant
and the $L^2$-norm of curvature (moreover, 
all of the the $L^2$-Betti numbers are bounded).
\end{remark}
 We emphasize that, in the Einstein case, the upper volume 
growth bound in (\ref{imp}) follows from the Bishop volume 
comparison theorem \cite{Bishop}, but since we are not assuming 
any pointwise Ricci curvature bound, this estimate is 
non-trivial. We state the volume comparison theorem
separately here, 
since it is of independent geometric interest:
\begin{theorem}
\label{bigthm_i}
Let $(X,g)$ be a complete, noncompact, $n$-dimensional 
Riemannian manifold with base point $p$.
Assume that there exists a constant $C_1 > 0$ so that
\begin{align}
\label{cond4_i}
Vol(B(q,s)) \geq C_1 s^n,
\end{align}
for any $q \in X$, and all $s \geq 0$.
Assume furthermore that as $r \rightarrow \infty$,
\begin{align}
\label{decay1_i}
\underset{S(r)}{sup} \ |Rm_g| &= o(r^{-2}),
\end{align}
where $S(r)$ denotes the sphere of radius $r$ centered at $p$. 
If $b_1(X) < \infty$, then $(X,g)$ has finitely many ends, 
and there exists a constant $C_2$ (depending on $g$) so that 
\begin{align}
\label{vga_i}
Vol(B(p,r)) \leq C_2 r^n.
\end{align}
Furthermore, each end is ALE of order $0$.
\end{theorem}

  An important problem is to find geometric conditions 
so that each end of a complete space will be
ALE of order $\tau > 0$, and to determine the 
optimal order of decay. We examine this problem for the following 
cases:\\
\begin{tabular}{ll}
\\
a. & Self-dual or anti-self-dual metrics with zero scalar curvature.\\

b. & Scalar-flat metrics with harmonic curvature.\\
\end{tabular}\\
\\
In class (a) we have scalar-flat 
K\"ahler metrics, which are an important class of extremal 
K\"ahler metrics \cite{Calabi1}. 
Class (b) is equivalent to scalar-flat and $\delta W = 0$, so 
in both classes we have scalar-flat locally conformally 
flat metrics \cite{SchoenYau1}, \cite{Schoen3}.
\begin{theorem}
\label{decayale}
Let $(X,g)$ be a complete, noncompact 4-dimensional Riemannian 
manifold with $g$ of class (a) or (b).
Assume that
\begin{align}
\int_X |Rm_g|^2 dV_g < \infty, \ C_S < \infty,
\mbox{ and }  b_1(X) < \infty.
\end{align}
Then $(X,g)$ has finitely many ends, and 
each end is ALE of order $\tau$ for any $\tau < 2$.
\end{theorem}

\begin{remark}
 It is unclear if this theorem 
can be improved exactly to ALE of order $2$, since this 
corresponds to an exceptional value, but 
we remark that ALE order of $2$ would be the best 
possible. When $|Rm| = O(r^{-4})$, 
in particular the mass is finite, but it is not necessarily zero. 
The examples of Claude LeBrun in \cite{LeBrun1}
are a particular case of (a), which have {\em{negative}} mass,
and  $|Rm| = O(r^{-4})$.
Another example is given by the 
Schwarzschild metric; on $\mathbf{R}^4 - \{0 \}$, let 
$g = (1 + m/(3r^2) )^{2} g_0.$
This is scalar-flat, and is by 
definition locally conformally flat, therefore 
it is Bach-flat, and also of harmonic curvature. For 
$m > 0$, the curvature decays like $C/r^4$, 
but the ADM mass is $m \neq 0$. 
Since $g$ has two ends, this example also shows that it is 
possible to have more than one end; this is 
in contrast to the Einstein case where the 
Cheeger-Gromoll splitting theorem rules out this
possibility. We remark also that in the Einstein case, it was proved in 
\cite{BKN} that $(X,g)$ is moreover ALE of order $4$.
\end{remark}

\begin{remark} An important application of the results 
in this paper is to the convergence of sequences of 
Bach-flat metrics on compact manifolds to orbifold metrics, 
this will be discussed in detail in a forthcoming paper. 
\end{remark} 

We end the introduction with a brief outline of the paper. 
In Section \ref{criticalequation}, we will 
discuss the systems satisfied by Bach-flat metrics 
and metrics with harmonic curvature. 
  The decay rate (\ref{imp2}) will be proved using a 
Moser iteration process in Section \ref{localregularity}, 
see Theorem \ref{higherlocalregthm2}. 
The crucial volume growth estimate in Theorem 
\ref{bigthm_i} will be proved in Section 
\ref{euclideanvolumegrowth}.

To prove Theorem \ref{decayale}, in Section \ref{proof2section} 
we will derive various improved Kato inequalities, 
which will allow us to improve 
the order of decay of the Ricci tensor to 
$|Ric| = O(r^{-\alpha})$ for any $\alpha < 4$.
In Section \ref{removal} we 
use a Yang-Mills type argument
as in \cite{Uhlenbeck2}, \cite{Tian}, to show 
that the full curvature tensor satisfies 
$|Rm| = O(r^{-\alpha})$ for any $\alpha < 4$. 
The result in \cite{BKN} then implies  
that $(X,g)$ is ALE of order $\tau$ for any 
$\tau < 2$. In Section~\ref{Constraints}, 
we remark that the Gauss-Bonnet and 
signature formulas put constraints on the 
ends of the ALE spaces which arise in Theorem \ref{decayale}. 

\subsection{Acknowledgements}

The authors are especially grateful to Tom Branson for 
improving our original Kato inequality in Section \ref{proof2section}
and showing us how the best constant follows easily from the 
general theory of Kato constants in \cite{Branson}. We would
also like to thank John Lott for assistance with the
developing map argument at the end of Section \ref{euclideanvolumegrowth}. The authors 
also benefitted on several occasions from conversations with
Toby Colding, Tom Mrowka, and Joao Santos.
We would also like to thank the referees for numerous 
suggestions which helped to improve the presentation.  

\section{Critical metrics}
\label{criticalequation}
In this section, we briefly discuss the systems of equations  
satisfied by Bach-flat metrics and metrics with harmonic 
curvature. 

\subsection{Bach-flat metrics with constant scalar curvature}
As stated in the introduction, the Euler-Lagrange equations 
of the functional
\begin{align}
\label{functional3}
\mathcal{W} : g \mapsto \int_X |W_g|^2 dV_g, 
\end{align}
in dimension $4$, are
\begin{align}
\label{bacheq2}
B_{ij} = \nabla^k \nabla^l W_{ikjl} +
\frac{1}{2}R^{kl}W_{ikjl} = 0.
\end{align}
Since the Bach tensor arises in the Euler-Lagrange 
equations of a Riemannian functional, it is symmetric,
and since the functional (\ref{functional3}) is 
conformally invariant, it follows that the Bach-flat 
equation (\ref{bacheq2}) is conformally invariant.

We note that (see \cite{Calderbankweaklyasd})
\begin{align}
B_{ij} = 2 \nabla^k \nabla^l W^+_{ikjl} +
R^{kl}W^+_{ikjl} = 2 \nabla^k \nabla^l W^-_{ikjl} +
R^{kl}W^-_{ikjl},
\end{align}
so that both self-dual and anti-self-dual metrics are 
Bach-flat. If $g$ is K\"ahler and has zero scalar 
curvature then $W^+ \equiv 0$ \cite{Derdzinski}, so 
these metrics are in particular anti-self-dual.  
Using the following Bianchi identities
\begin{align}
\label{Bianchi1}
\nabla^i R_{ijkl} = \nabla_k R_{jl} - \nabla_l R_{jk},
\end{align}
and 
\begin{align}
\label{Bianchi2}
\nabla^i W_{ijkl} = (n-3) (\nabla_k A_{jl} - \nabla_l A_{jk}), 
\end{align}
where $A_{ij}$ are the components of the Weyl-Schouten tensor 
\begin{align*}
A = \frac{1}{n-2} \left( Ric - \frac{1}{2(n-1)}R \cdot g \right),
\end{align*}
and $R$ denotes the scalar curvature, 
a computation shows that we may rewrite the 
Bach-flat equation (in dimension $4$) as
\begin{align*}
B_{ij} = \Delta A_{ij} - \nabla^k \nabla_i A_{jk}
+\frac{1}{2}R^{kl}W_{ikjl}= 0. 
\end{align*}
To simplify notation, in the following we will work 
in an orthonormal frame, and write all indices as lower.
We recall a standard formulas for commuting 
covariant derivatives: if $S$ is any 2-tensor, we have 
\begin{align}
\label{2tensor}
S_{ij;kl} - S_{ij;lk} = R_{klim}S_{mj} + R_{kljm}S_{im}.
\end{align}

Using (\ref{2tensor}), we have the formula 
\begin{align*}
\nabla_k \nabla_i A_{jk} = \nabla_i \nabla_k A_{jk}
 + R_{ikjp}A_{pk} + R_{ikkp}A_{jp}.
\end{align*}
So we may write the Bach tensor as
\begin{align*}
B_{ij} = \Delta A_{ij} - \nabla_i \nabla_k A_{jk}
+ R_{ikjp}A_{pk} - R_{ip}A_{jp}+\frac{1}{2}R_{kl}W_{ikjl}. 
\end{align*}
If $g$ also has constant scalar curvature, then 
\begin{align*}
\nabla_k A_{jk} = \frac{1}{6} \nabla_k R = 0, 
\end{align*}
so we have 
\begin{align*}
B_{ij} = \frac{1}{2}( \Delta Ric)_{ij} 
+ R_{ikjl}A_{kl} - R_{il}A_{jl} + A_{kl}W_{ikjl}.
\end{align*}
Therefore we may write the Bach-flat equations as
\begin{align}
(\Delta Ric)_{ij} = 2 ( R_{il}g_{jk} - R_{ikjl} -  W_{ikjl})A_{kl}. 
\end{align}
Introducing a convenient shorthand, we write this as 
\begin{align}
\label{shortcurveqn}
\Delta Ric = Rm* Ric.
\end{align}
\subsection{Metrics with harmonic curvature tensor}
The condition for harmonic curvature is 
that 
\begin{align}
\label{harmcond}
\delta Rm = -R_{ijkl;i} = 0. 
\end{align}
This condition was studied in \cite{Bo4}, 
\cite{Derd1}, \cite{Besse}, and is 
the Riemannian analogue of a Yang-Mills 
connection.  From the Bianchi identity (\ref{Bianchi1}), 
an equivalent condition is that $R_{ij;k} = R_{ik;j}$.
The Riemannian identity $R_{ik;i} = (1/2) R_{;k}$ implies
that the scalar curvature is constant.
An equivalent condition for harmonic curvature is therefore 
that $\delta W = 0$ and $ R = constant$. 
In particular, locally conformally flat 
metrics with constant scalar curvature have 
harmonic curvature. We compute
\begin{align*}
(\Delta Rm)_{ijkl} &= R_{ijkl;m;m}\\
&= ( - R_{ijlm;k} - R_{ijmk;l})_{;m}\\
&= -R_{ijlm;mk}  - R_{ijmk;ml} + Q(Rm)_{ijkl} = Q(Rm)_{ijkl}.
\end{align*}
where $Q(Rm)$ denotes a quadratic expression 
in the curvature tensor. 
In the shorthand, we write this as 
\begin{align}
\label{harmcurveqn}
\Delta Rm = Rm* Rm.
\end{align}
\section{Local regularity}
\label{localregularity}
We can consider more generally any system of the type
\begin{align}
\label{generaleqn}
\Delta Ric = Rm * Ric.
\end{align}
Any Riemannian metric satisfies 
\begin{align}
\label{curveqn}
\Delta Rm = L(\nabla^2 Ric) + Rm*Rm,
\end{align}
where $L(\nabla^2 Ric)$ denotes a linear expression 
in second derivatives of the Ricci tensor, and 
$Rm*Rm$ denotes a term which is quadratic in 
the curvature tensor (see \cite[Lemma 7.2]{Hamilton}). 

 For $X$ compact, we define the Sobolev constant $C_S$  
as the best constant $C_S$ so that 
for all $f \in C^{0,1}(X)$ (Lipschitz) we have
\begin{align}
\label{scc}
\Vert f \Vert_{L^{\frac{2n}{n-2}}} \leq C_S \left( 
\Vert \nabla f \Vert_{L^2} +  \Vert f \Vert_{L^2} \right).
\end{align}

  If $X$ is a complete, noncompact manifold, the Sobolev constant $C_S$ is 
defined as the best constant $C_S$ so that 
for all $f \in C^{0,1}_c(X)$ (Lipschitz with compact support),  we have
\begin{align}
\label{siq2}
\Vert f \Vert_{L^{\frac{2n}{n-2}}} \leq C_S \Vert \nabla f \Vert_{L^2}. 
\end{align}

 Even though second derivatives of the 
Ricci occur in (\ref{curveqn}), overall the principal
symbol of the system (\ref{generaleqn}) and (\ref{curveqn})
in triangular form. The equations (\ref{generaleqn}) and (\ref{curveqn}), 
when viewed as an elliptic system, together with the bound on the 
Sobolev constant, enable us to prove an $\epsilon$-regularity 
theorem (see \cite{Uhlenbeck1}, \cite{Uhlenbeck2}, \cite{Nakajima2},
\cite{Chang}):
\begin{theorem}
\label{higherlocalregthm}
Assume that (\ref{generaleqn}) is satisfied, 
let $r < diam(X)/2$, and $B(p,r)$ be a geodesic 
ball around the point $p$, and $k \geq 0$. Then there exist  
constants $\epsilon_0, C_k$  (depending upon $C_S$) so that if 
\begin{align*}
\Vert Rm \Vert_{L^2(B(p,r))} = 
\left\{ \int_{B(p,r)} |Rm|^2 dV_g \right\}^{1/2} \leq \epsilon_0,
\end{align*}
then 
\begin{align*}
\underset{B(p, r/2)}{sup}| \nabla^k Rm| \leq
\frac{C_k}{r^{2+k}} \left\{ \int_{B(p,r)} |Rm|^2 dV_g \right\}^{1/2}
\leq \frac{C_k \epsilon_0}{r^{2+k}}. 
\end{align*}
\end{theorem}
\begin{remark}
In the harmonic curvature case, the theorem
is much easier, since we have an equation on the full 
curvature tensor (\ref{harmcurveqn}), not just an equation on the Ricci;
see \cite{Akutagawa}, \cite{Anderson}, \cite{Nakajima2}, \cite{Tian}. 
\end{remark}
\begin{proof}
We will consider the case where the Sobolev inequality 
(\ref{scc}) is satisfied, the case of (\ref{siq2}) is similar. 
The proof will involve several lemmas. First, 
\begin{lemma}
\label{firstric}
There exist constants $\epsilon, C$ so that if 
$\Vert Rm \Vert_{L^2(B(p,r))} \leq \epsilon$, 
then 
\begin{align}
\label{needeqn}
\left\{ \int_{B(p,r/2)} |Ric|^{4} dV_g \right\}^{1/2}
\leq \frac{C}{r^2}\int_{B(p,r)} |Ric|^2 dV_g.
\end{align}
\end{lemma}
\begin{proof}
First,  from (\ref{generaleqn}), it follows that
\begin{align}
\label{Riceqn}
\Delta |Ric| \geq - |Rm| |Ric|. 
\end{align}
Without loss of generality, we may assume that $r=1$. 
The lemma then follows by scaling the metric. 
Let $0 \leq \phi \leq 1$ be a function supported in $B(p,1)$, then
\begin{align*}
\int_{B(p,1)} \phi^2 |Ric|^2 |Rm| &\geq 
\int_{B(p,1)} \phi^2 |Ric| ( - \Delta |Ric|)
= \int \nabla ( \phi^2 |Ric|) \cdot \nabla |Ric| \\
& = \int 2 \phi |Ric| \nabla \phi \cdot \nabla |Ric| 
+ \int \phi^2 |\nabla |Ric||^2 \\
& \geq  - \int (\delta^{-1}  | \nabla \phi| ^2 |Ric|^2 + 
{\delta}\phi^2 |\nabla|Ric||^2) +  \int | \phi \nabla |Ric||^2\\
& =   - \delta^{-1} \int | \nabla \phi| ^2 |Ric|^2 
+ (1 - \delta) \int | \phi \nabla |Ric||^2.
\end{align*}
Next, using the Sobolev constant bound, we have 
\begin{align*}
\left\{ \int ( \phi |Ric|)^4 \right\}^{1/2} 
&\leq C \int | \nabla ( \phi |Ric|)|^2
+ C \int \phi^2 |Ric|^2
\\
& = C \int | \nabla \phi|^2 |Ric|^2 + C \int | \phi \nabla |Ric||^2
+  C \int \phi^2 |Ric|^2.
\end{align*}
so choosing $\delta$ sufficiently small
\begin{align*}
\left\{ \int \phi^2 |Ric|^{4} \right\}^{1/2}
& \leq C \int \phi^2 |Ric|^2 |Rm| +  
C \int (\phi^2 +  | \nabla \phi| ^2) |Ric|^2 \\
&\leq C  \left\{
\int \phi^2 |Rm|^2 \right\}^{1/2}
\left\{
\int \phi^2 |Ric|^{4} \right\}^{1/2} 
+ C \int (\phi^2 + | \nabla \phi| ^2) |Ric|^2.
\end{align*}
Therefore for $\epsilon$ sufficiently small we have 
\begin{align}
\label{needeqnn}
\left\{ \int \phi^2 |Ric|^{4} \right\}^{1/2}
\leq  C \int ( \phi^2 + | \nabla \phi| ^2 ) |Ric|^2.
\end{align}
We then choose the cutoff function $\phi$ such that 
$\phi \equiv 1$ in $B(p,1/2)$, $\phi = 0 $ for 
$ r =1$, $|\nabla \phi| \leq C$, and we have 
\begin{align*}
\left\{ \int_{B(p,1/2)} |Ric|^{4} \right\}^{1/2}
\leq C\int_{B(p,1)} |Ric|^2.
\end{align*}
Scaling the metric, we obtain (\ref{needeqn}). 
\end{proof}
\begin{lemma}
There exist constants $\epsilon, C$ so that if 
$\Vert Rm \Vert_{L^2(B(p,r))} \leq \epsilon$,
then
\begin{align}
\label{mos1}
\int_{B(p,r/2)} |\nabla Ric|^2 dV_g 
\leq
\frac{C}{r^2} \int_{B(p,r)} |Ric|^2 dV_g,
 \end{align}
and
\begin{align}
\label{mos2}
\left\{ \int_{B(p,r/4)} |Rm|^4 dV_g \right\}^{1/2}
\leq
\frac{C}{r^2} \int_{B(p,r)} |Rm|^2 dV_g.
\end{align}
\end{lemma}
\begin{proof}
Again we may assume that $r=1$, and the lemma will follow 
by scaling the metric. 
Let $\phi$ be a cutoff function in $B(p,1)$, 
such that $\phi \equiv 1$ in $B(p, 1/2)$ and
$| \nabla \phi | \leq C$. We have 
\begin{align*}
&\int_{B(p,1)} \phi^2 |\nabla Ric|^2 = - \int \phi^2 \langle \Delta
Ric, Ric \rangle - 2 \int \phi \langle \nabla Ric, 
\nabla \phi \cdot Ric \rangle\\
& = - \int \phi^2 \langle Rm * Ric , Ric \rangle - 2 \int \phi
\langle \nabla Ric, \nabla \phi \cdot Ric \rangle\\
& \leq C \int \phi^2 |Rm| |Ric|^2 + \frac{C}{\delta} \int | \nabla
\phi|^2 |Ric|^2 + C \delta \int \phi^2 |\nabla Ric|^2\\
& \leq  C \left\{ \int \phi^2 |Rm|^2 \right\}^{1/2}
\left\{ \int \phi^2 |Ric|^4 \right\}^{1/2}
+  C \int_{B(p,1)} |Ric|^2 + C \delta \int \phi^2 |\nabla Ric|^2\\
& \leq  C \left\{ \int \phi^2 |Rm|^2 \right\}^{1/2}
\cdot C \int_{B(p,1)} |Ric|^2
+  C \int_{B(p,1)} |Ric|^2 + C \delta \int \phi^2 |\nabla Ric|^2,
\end{align*}
where we have used (\ref{needeqnn}).
By choosing $\delta$ small, and $\epsilon < 1$, we obtain
\begin{align*}
& \int_{B(p,1/2)} |\nabla Ric|^2 \leq 
(1 + \epsilon)  C \int_{B(p,1)} |Ric|^2
\leq 2 C \int_{B(p,1)} |Ric|^2.
\end{align*}
Scaling the metric, we obtain (\ref{mos1}). 
Next, let $\phi$ be a cutoff function in 
$B(p,1/2)$, such that $\phi \equiv 1$ in $B(p, 1/4)$ and
$| \nabla \phi | \leq C$, 
\begin{align*}
\int_{B(p,1/2)} & \langle \Delta Rm , \phi^2 Rm \rangle 
= \int  \langle \nabla^2 Ric + Rm* Rm, \phi^2 Rm \rangle \\
& =   \int  \langle \nabla^2 Ric, \phi^2 Rm \rangle 
+  \int \langle Rm* Rm, \phi^2 Rm \rangle \\
&=  - \int \langle 2 \phi \nabla Ric, \nabla \phi Rm \rangle
- \int  \phi^2 \langle \nabla Ric, \nabla Rm \rangle
+  \int \langle Rm* Rm, \phi^2 Rm \rangle.
\end{align*}
This yields 
\begin{align*}
\Big|  \int_{B(p,1/2)} & \langle \Delta Rm , \phi^2 Rm \rangle
\Big| \leq C \int   \phi^2 |\nabla Ric|^2  
+  C \int   |\nabla \phi|^2 |Rm|^2 \\
& + 
\frac{C}{\delta} \int  \phi^2 |\nabla Ric|^2
+ C \delta \int  \phi^2 |\nabla Rm|^2 +  C \int  \phi^2 |Rm| ^3.
\end{align*}
Integrating by parts,  
\begin{align*}
\int_{B(p,1/2)} \phi^2 | \nabla Rm|^2
& =  -\int \langle 2 \phi \nabla Rm, \nabla \phi \cdot Rm \rangle 
 - \int  \phi^2 \langle \Delta Rm , Rm \rangle\\
& \leq 
\frac{C}{\delta} \int  |\nabla \phi|^2 | Rm|^2
+ 2 C \delta \int  \phi^2 | \nabla Rm|^2   
+  C\int   |\nabla \phi|^2 |Rm|^2 \\
& \ \ \ \ \ + \frac{C'}{\delta} \int  \phi^2 |\nabla Ric|^2
+  C \int  \phi^2 |Rm| ^3.
\end{align*}
Choosing $\delta$ sufficiently small and using (\ref{mos1}), we obtain
\begin{align*}
\int \phi^2 |\nabla Rm|^2 
&\leq C \int | \nabla \phi|^2 |Rm|^2 + C \int_{B(p,1)} |Rm|^2
+ C \int \phi^2 |Rm|^3.
\end{align*}
Using the Sobolev inequality, 
\begin{align*}
\left\{ \int | \phi Rm|^4 \right\}^{1/2}
&\leq C \int | \nabla |\phi Rm| |^2 + C \int \phi^2 |Rm|^2\\
& \leq C \int | \nabla \phi|^2 |Rm|^2
+ C \int \phi^2 | \nabla |Rm||^2 + C \int \phi^2 |Rm|^2 \\
& \leq C \int ( \phi^2 + |\nabla \phi|^2) |Rm|^2
+ C \int \phi^2 | \nabla Rm|^2\\
& \leq  C \int_{B(p,1)} |Rm|^2 + C \int \phi^2 |Rm|^3\\
\leq C & \int_{B(p,1)} |Rm|^2
 + C \left\{ \int_{B(p,1/2)} |Rm|^2 \right\}^{1/2}
 \left\{ \int_{B(p,1/2)} \phi^4 |Rm|^4 \right\}^{1/2}. 
\end{align*}
Therefore by choosing $\epsilon$ small, we obtain 
\begin{align*}
\left\{ \int_{B(p,1/4)} |Rm|^4 \right\}^{1/2} \leq C 
\int_{B(p,1)} |Rm|^2.
\end{align*}
Scaling the metric, we obtain (\ref{mos2}). 
\end{proof}
\begin{lemma}
\label{iteration}
Suppose $0 \leq u \in L^2(B_{r})$ satisfies the inequality
\begin{align}
\Delta u \geq - u \cdot f
\end{align}
with 
\begin{align} 
\int_{B(p,r)} |f|^4 dV_g \leq \frac{C_1}{r^4}.
\end{align}
Then there exists a constant $C$ depending 
only upon $C_1, C_S$ so that  
\begin{align*}
\underset{B(p, r/2)}{sup} u \leq
\frac{C}{r^2} \Vert u \Vert_{L^2(B(p,r))}.
\end{align*}
\end{lemma}
\begin{proof}
This is a standard Moser iteration argument, 
see \cite[Lemma 4.6]{BKN}.
\end{proof}
\begin{lemma}
\label{supric}
There exist constants $\epsilon, C$ so that if 
$\Vert Rm \Vert_{L^2(B(p,r))} \leq \epsilon$, 
then 
\begin{align*}
\underset{B(p, r/2)}{sup}|Ric| \leq
\frac{C}{r^2} \Vert Ric \Vert_{L^2(B(p,r))}.
\end{align*}
\end{lemma}
\begin{proof}
From (\ref{generaleqn}), there exists a constant $C$ so that
\begin{align}
\Delta | Ric| \geq - C |Rm| |Ric|.
\end{align}
Also, using (\ref{mos2}), we have 
\begin{align}
\label{mos2'}
 \int_{B(p,r/2)} |Rm|^4 dV_g
\leq
\frac{C \epsilon}{r^4}
\end{align}
Therefore, we may apply Lemma \ref{iteration}, with 
$u = |Ric|$, $f = |Rm|$.
\end{proof}
An important consequence is the following
\begin{corollary}
\label{volbnd}
There exist constants $\epsilon, C$ so that if 
$\Vert Rm \Vert_{L^2(B(p,r))} \leq \epsilon$, 
then 
\begin{align*}
Vol(B(p,r/2)) \leq C r^4.
\end{align*}
\end{corollary}
\begin{proof}
Taking $\epsilon$ as in Lemma \ref{supric}, and scaling to 
unit size, we have 
\begin{align*}
\underset{B(p, 1/2)}{sup}|Ric| \leq
 C \cdot \epsilon.
\end{align*}
Therefore by the Bishop volume comparison 
theorem, $Vol(B(p,1/2)) \leq C$. 
\end{proof}
\begin{lemma}
There exist constants $\epsilon, C$ so that if 
$\Vert Rm \Vert_{L^2(B(p,r))} \leq \epsilon$,
then
\begin{align}
\label{moslater}
\int_{B(p,r/2)} |\nabla Ric|^2 dV_g 
\leq
\frac{C}{r^2} \int_{B(p,r)} |Ric|^2 dV_g.
\end{align}
\end{lemma}
\begin{remark}
This lemma will be crucial later in Section \ref{removal}. 
\end{remark}
\begin{proof}
Again we may assume that $r=1$, and the lemma will follow 
by scaling the metric. 
Let $\phi$ be a cutoff function in $B(p,1)$, 
such that $\phi \equiv 1$ in $B(p, 1/2)$ and
$| \nabla \phi | \leq C$. We have 
\begin{align*}
\int_{B(p,1)} & \phi^2 |\nabla Ric|^2
 = - \int \phi^2 \langle \Delta Ric,  Ric \rangle 
- 2 \int \langle \phi \nabla Ric,  \nabla \phi \cdot Ric \rangle\\
& =  - \int \phi^2 \langle Rm * Ric, Ric \rangle 
- 2 \int \langle \phi \nabla Ric , \nabla \phi \cdot Ric \rangle \\
& \leq C \int \phi^2 |Rm| |Ric|^2 + \frac{C}{\delta} \int | \nabla
\phi|^2 |Ric|^2 + C \delta \int \phi^2 |\nabla Ric|^2\\
& \leq  C \left\{ \int \phi^2 |Rm|^2 \right\}^{1/2}
\left\{ \int \phi^2 |Ric|^4 \right\}^{1/2}
+  C \int |Ric|^2 + C \delta \int \phi^2 |\nabla Ric|^2\\
& \leq  C \left\{ \int \phi^2 |Rm|^2 \right\}^{1/2}
\cdot C \int |Ric|^2
+  C \int |Ric|^2 + C \delta \int \phi^2 |\nabla Ric|^2,
\end{align*}
where we have used (\ref{needeqn}).
By choosing $\delta$ small, and $\epsilon < 1$, we obtain
\begin{align*}
& \int_{B(p,1/2)} |\nabla Ric|^2 \leq 
(1 + \epsilon)  C \int |Ric|^2
\leq 2 C \int |Ric|^2.
\end{align*}
\end{proof}
\begin{lemma}
\label{higherit}
For any $k \geq 0$, there exist constants $\epsilon, C$ so that if 
$\Vert Rm \Vert_{L^2(B(p,r))} \leq \epsilon$,
then 
\begin{align}
\label{lmos1}
\left\{ \int_{B(p,r/2)} |\nabla^{k} Ric|^{4} dV_g \right\}^{1/2}
\leq \frac{C}{r^{2(k+1)}} \int_{B(p,r)} |Rm|^2 dV_g, 
\end{align}
\begin{align}
\label{lmos2}
\int_{B(p,r/4)} |\nabla^{k+1} Ric|^2 dV_g
\leq \frac{C}{r^{2(k+1)}} \int_{B(p,r)} |Rm|^2 dV_g,
\end{align} 
\begin{align}
\label{lmos3}
\int_{B(p,r/8)} |\nabla^k Rm|^2 dV_g
\leq \frac{C}{r^{2k}} \int_{B(p,r)} |Rm|^2 dV_g,
\end{align}
\begin{align}
\label{lmos4}
\left\{ \int_{B(p,r/16)} |\nabla^k Rm|^4 dV_g \right\}^{1/2}
\leq \frac{C}{r^{2(k+1)}} \int_{B(p,r)} |Rm|^2 dV_g.
\end{align}
\end{lemma}
\begin{proof} 
We use induction on $k$. The case $k=0$ has been proved in 
the previous lemmas. Assume the lemma is true 
for $i = 1 \dots k-1$.
In the following proof, one must shrink 
the ball further after each step of the iteration.
For simplicity, we will do this automatically without
mention, so that the final step shrinks by the 
appropriate factor.

 From the equation $ \Delta Ric = Rm * Ric$, it follows that
(\cite[Theorem 2.4]{Chang}) 
\begin{align}
\label{deltaRic}
\Delta ( \nabla^k Ric) = \sum_{l=0}^{k} \nabla^l Rm * \nabla^{k-l} Ric.
\end{align}
Similarly, from the equation $ \Delta Rm = Rm* Rm + \nabla^2
Ric$, it follows that (\cite[Corollary 3.11]{Chang}) 
\begin{align}
\label{deltaRm}
\Delta ( \nabla^k Rm) = \sum_{l=0}^{k} \nabla^l Rm * \nabla^{k-l}
Rm + \nabla^{k+2} Ric. 
\end{align}
\noindent{\em{Proof of (\ref{lmos1})}}.
We let $\phi$ be a cutoff function as before, and consider 
the following expression
\begin{align}
\notag
\int_{B(p,1)}& \langle \Delta( \nabla^k Ric) , \phi^2
 \nabla^k Ric \rangle 
 = \int  \langle  \sum_{l=0}^{k} \nabla^l Rm * \nabla^{k-l} Ric  , \phi^2
 \nabla^k Ric \rangle\\
\begin{split}
\label{gg1}
&= \int  \langle \nabla^k Rm * Ric  , \phi^2 \nabla^k Ric \rangle
+  \int  \langle Rm * \nabla^k Ric  , \phi^2 \nabla^k Ric \rangle\\
&\ \ \ \ + \sum_{l=1}^{k-1} \int  \langle \nabla^l Rm * \nabla^{k-l} Ric  , \phi^2
 \nabla^k Ric \rangle.
\end{split}
\end{align}
Integrate the first term in (\ref{gg1}) by parts
\begin{align*}
& \int  \langle \nabla^k Rm * Ric  , \phi^2 \nabla^k Ric \rangle
=  -2 \int  \langle  \nabla \phi \cdot \nabla^{k-1} Rm * Ric, 
\phi \nabla^k Ric \rangle\\
&- \int  \langle  \nabla^{k-1} Rm * \nabla Ric, \phi^2
 \nabla^k Ric \rangle 
-  \int  \langle  \nabla^{k-1} Rm * Ric, \phi^2
 \nabla^{k+1} Ric \rangle.
\end{align*}
We estimate
\begin{align*} 
&\Big|\int  \langle \nabla^k Rm * Ric  , \phi^2 \nabla^k Ric \rangle \Big|\\
&\leq C \int \phi |\nabla \phi| | \nabla^{k-1} Rm| | Ric| |\nabla^k Ric|
+C  \int \phi^2 | \nabla^{k-1} Rm| |\nabla Ric| |\nabla^k Ric| \\
& \ \ \ + C \int  \phi^2 |\nabla^{k-1} Rm| |Ric| |\nabla^{k+1} Ric|\\
&\leq  C \int \phi^2 | \nabla^{k-1} Rm|^4 
+ C \int \phi^2 | Ric|^4 
+ C \int \phi^2 |\nabla^k Ric|^2\\
& \ \ \ \ +C  \int \phi^2 | \nabla^{k-1} Rm| |\nabla Ric| |\nabla^k Ric|  \\
& \ \ \ \ +  \frac{C}{\delta} \Big( \int  \phi^2 |\nabla^{k-1} Rm|^4 
+ \int \phi^2 |Ric|^4 \Big) + C \delta \int  \phi^2 |\nabla^{k+1} Ric|^2.
\end{align*}
Using the inductive hypothesis, we obtain
\begin{align*}\Big|\int  &\langle \nabla^k Rm * Ric  , \phi^2
  \nabla^k Ric \rangle \Big| \\
 &\leq C \delta \int  \phi^2 |\nabla^{k+1} Ric|^2
+ C \int | Rm|^2 + C \Big\{ \int |Rm|^2 \Big\}^2\\
&+C  \int \phi^2 | \nabla^{k-1} Rm| |\nabla Ric| |\nabla^k Ric| 
\end{align*}
Therefore when $\Vert Rm \Vert_{L_2(B(p,1))} < \epsilon < 1$, 
we have
\begin{align*}
\Big|\int  &\langle \nabla^k Rm * Ric  , \phi^2
  \nabla^k Ric \rangle \Big| 
\leq C \delta \int  \phi^2 |\nabla^{k+1} Ric|^2
+ C \int |Rm|^2 \\
&+C  \int \phi^2 | \nabla^{k-1} Rm| |\nabla Ric| |\nabla^k Ric| .
\end{align*}
Next we estimate the second term in (\ref{gg1})
\begin{align*} 
\int  & \langle Rm * \nabla^k Ric  , \phi^2 \nabla^k Ric \rangle
\leq  \Big\{ \int \phi^2 |Rm|^2 \Big\}^{1/2}
\Big\{ \int \phi^2 |\nabla^k Ric|^4 \Big\}^{1/2}\\
& \leq \epsilon \Big\{ \int \phi^2 |\nabla^k Ric|^4 \Big\}^{1/2}
\end{align*}
Next we estimate the sum term in (\ref{gg1}) (for $k \geq 2$)
\begin{align*}
&\sum_{l=1}^{k-1} \int  \langle \nabla^l Rm * \nabla^{k-l} Ric  , \phi^2
 \nabla^k Ric \rangle \\
&\leq \sum_{l=1}^{k-1} \Big( \int \phi^2 | \nabla^l Rm |^4 
+ \int \phi^2 |\nabla^{k-l} Ric|^4 + \int \phi^2 |\nabla^k Ric|^2 \Big)
\leq C \int |Rm|^2.
\end{align*}
by the inductive hypothesis, and since $\Vert Rm \Vert_{L^2(B(p,1))} <
\epsilon < 1$. 
Adding these estimates we obtain
\begin{align*}
\Big| \int_{B(p,1)}& \langle \Delta( \nabla^k Ric) , \phi^2
 \nabla^k Ric \rangle \Big | \\
 &\leq C \delta \int  \phi^2 |\nabla^{k+1} Ric|^2
+ C \int | Rm|^2 \\
&+ C  \int \phi^2 | \nabla^{k-1} Rm| |\nabla Ric| |\nabla^k Ric|
+  \epsilon \Big\{ \int \phi^2 |\nabla^k Ric|^4 \Big\}^{1/2}.
\end{align*}
On the other hand, we have 
\begin{align*}
\int_{B(p,1)} \langle \Delta \nabla^k Ric,  \phi^2 \nabla^k Ric 
\rangle 
&= - 2 \int \langle \phi \nabla^{k+1} Ric, \nabla \phi \cdot \nabla^k Ric 
\rangle \\
& \ \ \ \ \ \ - \int \langle \nabla^{k+1} Ric, \phi^2 \nabla^{k+1} Ric \rangle,
\end{align*}
which implies that 
\begin{align*}
& \int \phi^2 |\nabla^{k+1} Ric|^2
 \leq \Big| \int \langle \Delta \nabla^k Ric , \phi^2
\nabla^k Ric \rangle \Big| 
+ 2 \int \phi | \nabla \phi| | \nabla^k Ric| |\nabla^{k+1} Ric|\\
& \leq \Big| \int \langle \Delta \nabla^k Ric , \phi^2
\nabla^k Ric \rangle \Big| 
+ \frac{2}{\delta} \int | \nabla \phi|^2 | \nabla^k Ric|^2 
+ 2 \delta \int \phi^2 |\nabla^{k+1} Ric|^2\\
& \leq  C \delta \int  \phi^2 |\nabla^{k+1} Ric|^2
+ C \int | Rm|^2 \\
&+ C  \int \phi^2 | \nabla^{k-1} Rm| |\nabla Ric| |\nabla^k Ric|
+  \epsilon \Big\{ \int \phi^2 |\nabla^k Ric|^4 \Big\}^{1/2}\\
& + \frac{2}{\delta} \int | \nabla \phi|^2 | \nabla^k Ric|^2 
+ 2 \delta \int \phi^2 |\nabla^{k+1} Ric|^2.
\end{align*}
So by choosing $\delta$ small, 
\begin{align*}
& \int \phi^2 |\nabla^{k+1} Ric|^2
\leq C \int | Rm|^2 + 
C  \int \phi^2 | \nabla^{k-1} Rm| |\nabla Ric| |\nabla^k Ric|\\
& +  \epsilon \Big\{ \int \phi^2 |\nabla^k Ric|^4 \Big\}^{1/2}
+ C \int | \nabla \phi|^2 | \nabla^k Ric|^2. 
\end{align*}
Next from the Sobolev inequality and Lemma (\ref{supric}), 
\begin{align*}
\left\{ \int ( \phi | \nabla^{k} Ric|)^4 \right\}^{1/2}
&\leq C \int | \nabla ( \phi |\nabla^{k} Ric|)|^2 + C\int \phi^2
| \nabla^{k} Ric |^2  \\
& \leq C \int ( \phi^2 + | \nabla \phi|^2) | \nabla^{k} Ric|^{2} 
+ C \int \phi^2 |\nabla^{k+1} Ric|^2 \\
& \leq C \int ( \phi^2 + | \nabla \phi|^2 ) | \nabla^k Ric|^{2} 
+  C \int | Rm|^2 \\
& + C  \int \phi^2 | \nabla^{k-1} Rm| |\nabla Ric| |\nabla^k Ric|
 +  \epsilon \Big\{ \int \phi^2 |\nabla^k Ric|^4 \Big\}^{1/2}. 
\end{align*}
Choosing $\epsilon$ small, and using the inductive
hypothesis, we obtain
\begin{align*}
\left\{ \int ( \phi | \nabla^{k} Ric|)^4 \right\}^{1/2}
& \leq C \int | Rm|^2 +
 + C  \int \phi^2 | \nabla^{k-1} Rm| |\nabla Ric| |\nabla^k Ric|. 
\end{align*}
Now if $k=1$, the last term is 
\begin{align*}
\int \phi^2 | Rm| |\nabla Ric|^2
\leq \left\{ \int_{B(p,1)} | Rm|^2 \right\}^{1/2}
\left\{ \int \phi^4 |\nabla Ric|^4 \right \}^{1/2}
\leq \epsilon \left\{ \int \phi^4 |\nabla Ric|^4 \right \}^{1/2},
\end{align*}
so this term may be absorbed into the left, and we obtain
\begin{align}
\label{gradmos}
\left\{ \int ( \phi | \nabla Ric|)^4 \right\}^{1/2}
& \leq C \int | Rm|^2.
\end{align}
For $k \geq 2$, the last term is 
\begin{align*}
\int \phi^2 & | \nabla^{k-1} Rm| |\nabla Ric| |\nabla^k Ric|
\leq C\int \phi^2 | \nabla^{k-1} Rm|^2 |\nabla Ric|^2
+ C \int \phi^2 |\nabla^k Ric|^2 \\
& \leq C \int \phi^2 | \nabla^{k-1} Rm|^4 
+ C \int \phi^2 |\nabla Ric|^4
+ C \int \phi^2 |\nabla^k Ric|^2,
\end{align*}
so by the inductive hypothesis and $\Vert Rm \Vert_{L^2(B(p,1))} < \epsilon
<1$, we obtain
\begin{align}
\label{gradmos0}
\left\{ \int ( \phi | \nabla^k Ric|)^4 \right\}^{1/2}
& \leq C \int | Rm|^2.
\end{align}
Inequality (\ref{lmos1}) follows by choosing $\phi$ similar
to before and then scaling the metric. 

\noindent {\em{Proof of (\ref{lmos2}).}}
Integrating by parts and using (\ref{deltaRic}), 
\begin{align}
\notag
 \int_{B(p,r)} & \phi^2 | \nabla^{k+1} Ric|^2 
= - \int \langle \phi^2 \nabla^k Ric , \Delta \nabla^k Ric \rangle 
- 2 \int \langle \nabla \phi \nabla^k Ric  , \phi \nabla^{k+1} Ric
\rangle\\
\label{gg2}
& =   - \int \langle \phi^2 \nabla^k Ric , 
\sum_{l=0}^{k} \nabla^l Rm * \nabla^{k-l} Ric \rangle  
- 2 \int \langle \nabla \phi \nabla^k Ric  , \phi \nabla^{k+1} Ric \rangle.
\end{align}
For the first term in (\ref{gg2}), 
\begin{align*}
& -\sum_{l=0}^{k} \int \langle \phi^2 \nabla^k Ric , 
 \nabla^l Rm * \nabla^{k-l} Ric \rangle \\
& = \int \langle \phi^2 \nabla^k Ric , 
 \nabla^k Rm * Ric \rangle 
-\sum_{l=0}^{k-1} \int \langle \phi^2 \nabla^k Ric , 
 \nabla^l Rm * \nabla^{k-l} Ric \rangle. 
\end{align*}
Integrating by parts 
\begin{align*}
& \int \langle \phi^2 \nabla^k Ric , \nabla^k Rm * Ric \rangle
 = - 2 \int \langle \phi \nabla^k Ric , \nabla \phi \cdot \nabla^{k-1} Rm 
* Ric \rangle\\
& - \int \langle \phi^2 \nabla^{k+1} Ric , \nabla^{k-1} Rm * Ric \rangle
- \int \langle \phi^2 \nabla^k Ric , \nabla^{k-1} Rm * \nabla Ric \rangle,
\end{align*}
consequently 
\begin{align*}
& \Big| \int \langle \phi^2 \nabla^k Ric , \nabla^k Rm * Ric \rangle \Big|
\leq C \int |\phi| |\nabla \phi| |\nabla^k Ric| |\nabla^{k-1} Rm| 
|Ric| \\
& + C \int \phi^2 |\nabla^{k+1} Ric| |\nabla^{k-1} Rm| | Ric |
+ C\int \phi^2 |\nabla^k Ric| |\nabla^{k-1} Rm|
|\nabla Ric| \\
& \leq \int \phi^2 |\nabla^k Ric|^2 
+ \int |\nabla \phi|^2 |\nabla^{k-1} Rm|^2 |Ric|^2 \\
&\ \ \ \ + C \delta  \int \phi^2 |\nabla^{k+1} Ric|^2 
+ \frac{C}{\delta} \int \phi^2 |\nabla^{k-1} Rm|^2 |Ric|^2\\
&\ \ \ \  + C \int \phi^2 |\nabla^k Ric|^2 
+ C\int \phi^2 |\nabla^{k-1} Rm|^2 |\nabla Ric|^2\\
& \leq C \int |Rm|^2  +  C \delta  \int \phi^2 |\nabla^{k+1} Ric|^2, 
\end{align*}
by the inductive hypothesis and (\ref{lmos1}).
Also by induction, we have 
\begin{align*}
& \Big| \sum_{l=0}^{k-1} \int \langle \phi^2 \nabla^k Ric ,  \nabla^l
 Rm * \nabla^{k-l} Ric \rangle \Big|\\
& \leq \sum_{l=0}^{k-1} \Big( \int | \nabla^l Rm |^4 
+ \int |\nabla^{k-l} Ric|^4 + \int |\nabla^k Ric|^2 \Big)
\leq C \int |Rm|^2.
\end{align*}
We estimate the last term in (\ref{gg2}) 
\begin{align*}
&\Big| 2 \int \langle \nabla \phi \nabla^k Ric  , \phi \nabla^{k+1} Ric
\rangle \Big|
\leq C \int |\nabla \phi| \phi |\nabla^k Ric| | \nabla^{k+1} Ric|\\
&\leq C \delta \int \phi^2 | \nabla^{k+1} Ric|^2
+ \frac{C}{\delta} 
\int |\nabla \phi|^2 |\nabla^k Ric|^2\\
&  \leq  C \delta \int \phi^2 | \nabla^{k+1} Ric|^2+ C \int |Rm|^2.
\end{align*}
Combining the above estimates, we obtain
\begin{align*}
 \int_{B(p,r)} & \phi^2 | \nabla^{k+1} Ric|^2 
\leq  C \delta \int \phi^2 | \nabla^{k+1} Ric|^2
+  C \int |Rm|^2,
\end{align*}
So by choosing $\delta$ small, and $\phi$ as before, we obtain (\ref{lmos2}).

\noindent{\em{Proof of (\ref{lmos3}).}}
Integrating by parts and using (\ref{deltaRm}), 
\begin{align*}
 \int_{B(p,1)} \phi^2 | \nabla^{k} Rm|^2 
& = - \int \langle \phi^2 \nabla^{k-1} Rm , 
\Delta \nabla^{k-1} Rm \rangle 
- 2 \int \langle \nabla \phi \nabla^{k-1} Rm  , \phi \nabla^{k} Rm
\rangle\\
& =   - \int \langle \phi^2 \nabla^{k-1} Rm , 
\sum_{l=0}^{k-1} \nabla^l Rm * \nabla^{k-1-l} Rm + \nabla^{k+1} Ric 
\rangle\\
& \ \ \ \ - 2 \int \langle \nabla \phi \nabla^{k-1} Rm, 
\phi \nabla^{k} Rm \rangle,
\end{align*}
which implies the estimate
\begin{align*}
& \int_{B(p,1)} \phi^2 | \nabla^{k} Rm|^2 
\leq   \sum_{l=0}^{k-1} \int \phi^2 |\nabla^{k-1} Rm|
|\nabla^l Rm| |\nabla^{k-1-l} Rm| \\
& + \int \phi^2 |\nabla^{k-1} Rm| |\nabla^{k+1} Ric| 
+  2 \int \phi |\nabla \phi| |\nabla^{k-1} Rm| |\nabla^{k} Rm|\\
&\leq C \sum_{l=0}^{k-1} \Big( \int \phi^2 |\nabla^{k-1} Rm|^4
+ \int \phi^2 |\nabla^l Rm|^4 + \int \phi^2 |\nabla^{k-1-l} Rm|^2 \Big)\\ 
& + \left\{ \int \phi^2 |\nabla^{k-1} Rm|^2 | \right\}^{1/2} 
 \left\{ \int \phi^2 |\nabla^{k+1} Ric|^2 \right\}^{1/2}\\
& + \frac{C}{\delta}  \int |\nabla \phi|^2 |\nabla^{k-1} Rm|^2  
+ C \delta  \int \phi^2 |\nabla^{k} Rm|^2.
\end{align*}
Using the inductive hypothesis, (\ref{lmos2}), 
and choosing $\delta$ sufficiently
small, we obtain
\begin{align*}
 \int_{B(p,1)} \phi^2 | \nabla^{k} Rm|^2 
\leq C \int |Rm|^2 .
\end{align*}

\noindent{\em{Proof of (\ref{lmos4})}}.
We let $\phi$ be a cutoff function as before, and consider 
the following expression
\begin{align}
\notag
\int_{B(p,1)}& \langle \Delta( \nabla^k Rm) , \phi^2
 \nabla^k Rm \rangle \\
\notag
& = \sum_{l=0}^{k} \int  \langle \nabla^l Rm * \nabla^{k-l} Rm + 
\nabla^{k+2} Ric  , \phi^2 \nabla^k Rm \rangle\\
\begin{split}
\label{gg4}
&= \int  \langle \nabla^k Rm * Rm  , \phi^2 \nabla^k Rm \rangle
+  \int  \langle  \nabla^{k+2} Ric  , \phi^2 \nabla^k Rm \rangle\\
&+ \sum_{l=1}^{k-1} \int  \langle \nabla^l Rm * \nabla^{k-l} Rm  , \phi^2
 \nabla^k Rm \rangle.
\end{split}
\end{align}
For the first term in (\ref{gg4}),
\begin{align*}
& \Big| \int  \langle \nabla^k Rm * Rm  , \phi^2 \nabla^k Rm
 \rangle \Big|
\leq C \int \phi^2 |Rm| |\nabla^k Rm|^2\\
&\leq  \left\{ \int \phi^2 |Rm^2| \right\}^{1/2}   
\left\{ \int \phi^2 |\nabla^k Rm|^4 \right\}^{1/2}
\leq \epsilon \left\{ \int \phi^2 |\nabla^k Rm|^4 \right\}^{1/2}.
\end{align*}
Next we estimate the second term in (\ref{gg4}) 
\begin{align*} 
\int  \langle \nabla^{k+2}  Ric  , \phi^2 \nabla^k Rm \rangle
 =- & \int   \langle \nabla^{k+1}  Ric  , \phi^2 \nabla^{k+1} Rm \rangle\\
& - \int  \langle \nabla^{k+1}  Ric  , 2\phi \nabla \phi \nabla^k 
Rm \rangle.
\end{align*}
So we have 
\begin{align*} 
\Big| \int   \langle \nabla^{k+2}  Ric  , \phi^2 \nabla^k Rm \rangle \Big|
& \leq C \int  |\nabla^{k+1}  Ric | \phi^2 |\nabla^{k+1} Rm| \\
& + C \int   |\nabla^{k+1}  Ric| \phi |\nabla \phi| |\nabla^k Rm|\\
&\leq \frac{C}{\delta} \int  \phi^2 |\nabla^{k+1}  Ric |^2
+ C \delta \int \phi^2 |\nabla^{k+1} Rm|^2\\
& + C \int  \phi^2 |\nabla^{k+1}  Ric|^2 
+ C  \int |\nabla \phi|^2 |\nabla^k Rm|^2.
\end{align*}
Using induction and (\ref{lmos2}) we obtain
\begin{align*} 
\Big| \int   \langle \nabla^{k+2}  Ric  , \phi^2 \nabla^k Rm \rangle \Big|
& \leq C \delta \int \phi^2 |\nabla^{k+1} Rm|^2 +  C  \int |Rm|^2.
\end{align*}
Next we estimate the sum term in (\ref{gg4}) (for $k \geq 2$)
\begin{align*}
&\sum_{l=1}^{k-1} \int  \langle \nabla^l Rm * \nabla^{k-l} Rm  , \phi^2
 \nabla^k Rm \rangle \\
&\leq \sum_{l=1}^{k-1} \Big( \int \phi^2 | \nabla^l Rm |^4 
+ \int \phi^2 |\nabla^{k-l} Rm|^4 + \int \phi^2 |\nabla^k Rm|^2 \Big)
\leq C \int |Rm|^2.
\end{align*}
by (\ref{lmos3}), the inductive hypothesis, and since $\Vert Rm \Vert_{L^2(B(p,1))} <
\epsilon < 1$. 
Adding these estimate we obtain
\begin{align*}
\Big| \int_{B(p,1)}& \langle \Delta( \nabla^k Rm) , \phi^2
 \nabla^k Rm \rangle \Big | \\
 &\leq C \delta \int  \phi^2 |\nabla^{k+1} Rm|^2
+ C \int | Rm|^2 + \epsilon \left\{ \int \phi^2 |\nabla^k Rm|^4 \right\}^{1/2}.
\end{align*}
On the other hand, we have 
\begin{align*}
\int_{B(p,1)} \langle \Delta \nabla^k Rm,  \phi^2 \nabla^k Rm 
\rangle 
= &- \int \langle \nabla^{k+1} Rm, 2 \phi \nabla \phi \nabla^k Rm 
\rangle\\  
&\ \ \ \ \ - \int \langle \nabla^{k+1} Rm, \phi^2 \nabla^{k+1} Rm \rangle,
\end{align*}
which implies that 
\begin{align*}
& \int \phi^2 |\nabla^{k+1} Rm|^2
 \leq \Big| \int \langle \Delta \nabla^k Rm , \phi^2
\nabla^k Rm \rangle \Big| 
+ 2 \int \phi | \nabla \phi| | \nabla^k Rm| |\nabla^{k+1} Rm|\\
& \leq \Big| \int \langle \Delta \nabla^k Rm , \phi^2
\nabla^k Rm \rangle \Big| 
+ \frac{2}{\delta} \int | \nabla \phi|^2 | \nabla^k Rm|^2 
+ 2 \delta \int \phi^2 |\nabla^{k+1} Rm|^2\\
& \leq  C \delta \int  \phi^2 |\nabla^{k+1} Rm|^2
+ C \int | Rm|^2 +  \epsilon \left\{ \int \phi^2 |\nabla^k Rm|^4 \right\}^{1/2} \\
& \ \ \ \ \ + \frac{2}{\delta} \int | \nabla \phi|^2 | \nabla^k Rm|^2 
+ 2 \delta \int \phi^2 |\nabla^{k+1} Rm|^2.
\end{align*}
So by choosing $\delta$ small, we have 
\begin{align*}
& \int \phi^2 |\nabla^{k+1} Rm|^2
\leq C \int | Rm|^2 +   \epsilon \left\{ \int \phi^2 |\nabla^k Rm|^4
\right\}^{1/2}
+ C \int | \nabla \phi|^2 | \nabla^k Rm|^2. 
\end{align*}
Next from the Sobolev inequality we have 
\begin{align*}
& \left\{ \int ( \phi | \nabla^{k} Rm|)^4 \right\}^{1/2}
\leq C \int | \nabla ( \phi |\nabla^{k} Rm|)|^2
+ C \int \phi^2 | \nabla^{k} Rm|^2 \\
& \leq C \int ( \phi^2 + | \nabla \phi|^2) | \nabla^{k} Rm|^{2} 
+ C \int \phi^2 |\nabla^{k+1} Rm|^2 \\
& \leq C \int (\phi^2 + | \nabla \phi|^2) | \nabla^k Rm|^{2} 
+  C \int | Rm|^2 + \epsilon \left\{ \int \phi^2 |\nabla^k Rm|^4 \right\}^{1/2}. 
\end{align*}
Using the inductive hypothesis and choosing $\epsilon$ sufficiently 
small, we obtain
\begin{align}
\label{gradmos0non}
\left\{ \int ( \phi | \nabla^k Rm|)^4 \right\}^{1/2}
& \leq C \int | Rm|^2.
\end{align}
Inequality (\ref{lmos4}) follows by scaling the metric. 

\end{proof}
The following is the main Moser iteration lemma
\begin{lemma}
\label{iteration1}
Suppose $0 \leq u \in L^2(B(p,r))$ satisfies the inequality
\begin{align}
\Delta u \geq - u \cdot f - h,
\end{align}
with $h \geq 0$ and $f \geq 0$ satisfying 
\begin{align} 
\int_{B(p,r)} f^4 dV_g \leq \frac{C_1}{r^4}.
\end{align}
If $Vol(B(p,r)) \leq C_2 r^4$, 
then there exists a constant $C$ depending 
only upon $C_1, C_2, C_S$ so that 
\begin{align*}
\underset{B(p, r/2)}{sup} u \leq
\frac{C}{r^2} \Vert u \Vert_{L^2(B(p,r))}
+ C r \Vert h \Vert_{L^4(B(p,r))}.
\end{align*}
\end{lemma}
\begin{proof}
Let $k = r \Vert h \Vert_{L^4(B(p,r))}$,
and $\bar{u} = u + k$. 
Assume $k \neq 0$, then $\bar{u}$ satisfies 
\begin{align*}
\Delta \bar{u} = \Delta u \geq -uf - h \geq -\bar{u} f -h 
\geq - \bar{u} \Big( f + \frac{h}{k} \Big)= - \bar{u}\bar{f}.
\end{align*}
We have 
\begin{align*} 
\left\{ \int_{B(p,r)} \bar{f}^4 dV_g \right\}^{1/4}
&\leq  \left\{ \int_{B(p,r)} f^4 dV_g \right\}^{1/4}
+ \left\{ \int_{B(p,r)} \Big( \frac{h}{k}\Big)^4 dV_g \right\}^{1/4}\\
&\leq \frac{C_1}{r}
+ \frac{1}{k}\left\{ \int_{B(p,r)} h^4 dV_g \right\}^{1/4}\\
&\leq \frac{(C_1+1)}{r}.
\end{align*}
Lemma \ref{iteration} implies that 
\begin{align*}
\underset{B(p, r/2)}{sup} \bar{u} \leq
\frac{C}{r^2} \Vert \bar{u} \Vert_{L^2(B(p,r))}.
\end{align*}
In terms of $u$, we have 
\begin{align*}
\underset{B(p, r/2)}{sup} u &\leq
\frac{C}{r^2} \Vert \bar{u} \Vert_{L^2(B(p,r))}
\leq \frac{C}{r^2} \Big(  \Vert u \Vert_{L^2(B(p,r))}
+ \Vert k \Vert_{L^2(B(p,r))} \Big)\\
& \leq \frac{C}{r^2} \Big(  \Vert u \Vert_{L^2(B(p,r))}
+ k (Vol(B(p,r)))^{1/2} \Big) \\
& \leq \frac{C}{r^2} \Vert u \Vert_{L^2(B(p,r))}
+ Cr \Vert h \Vert_{L^4(B(p,r))}. 
\end{align*}
\end{proof}
We now complete the proof of the Theorem. We may assume 
that $r=1$. First consider the 
case $k=0$.
From (\ref{curveqn}), it follows there 
exists constants $C_1$ and $C_2$ so that 
\begin{align}
\label{Rmeqn}
\Delta |Rm| \geq -C_1 |Rm|^2 - C_2 | \nabla^2 Ric|.
\end{align}
This case $k=0$ follows by applying Lemma \ref{iteration1} to
(\ref{Rmeqn}),
using the $L^4$ bound on $\nabla^2 Ric$ from Lemma \ref{higherit},
and the volume bound from Corollary \ref{volbnd}. 

 Assume the theorem is true for $i = 0 \dots k-1$, 
so that we have pointwise bounds on $| \nabla^i Rm|$ 
for $ i = 0 \dots k-1$. 
From Lemma \ref{higherit}, we have 
bounds on $\Vert \nabla^{k+2} Ric \Vert_{L^{4}}$, 
Using (\ref{deltaRm}), we obtain
\begin{align*}
\Delta| \nabla^k Rm| &\geq - C \sum_{l=0}^{k} |\nabla^l Rm| 
|\nabla^{k-l}Rm| -  C |\nabla^{k+2} Ric |\\
&=  -C | \nabla^k Rm| |Rm| - C \sum_{l=1}^{k-1} |\nabla^l Rm| 
|\nabla^{k-l}Rm| -  C |\nabla^{k+2} Ric |.
\end{align*}
Lemma \ref{iteration1} and Corollary \ref{volbnd} 
then yield the pointwise bounds on $|\nabla^k Rm|$.
\end{proof}

We next apply Theorem \ref{higherlocalregthm} to 
noncompact spaces to give a rate of curvature decay at
infinity. 
\begin{theorem}
\label{higherlocalregthm2} 
Assume that $(X,g)$ is a complete, noncompact space 
satisfying (\ref{generaleqn}), 
\begin{align}
\label{a1'}
\int_X |Rm_g|^2 dV_g < \infty, \mbox{ and } C_S < \infty.
\end{align}
  Fix a basepoint $p$, and let 
$D(r) = X \setminus B(p,r)$ be the complement of a geodesic 
ball around the point $p$, and $k \geq 0$. There exist  
constants $\epsilon_0, C_k$  (depending upon $C_S$) so that if 
$\Vert Rm \Vert_{L^2(D(r))} \leq \epsilon_0,$
then 
\begin{align*}
\underset{D(2r)}{sup}| \nabla^k Rm| \leq
\frac{C_k}{r^{2+k}} \left\{ \int_{D(r)} |Rm|^2 dV_g \right\}^{1/2}
\leq \frac{C_k \epsilon_0}{r^{2+k}}. 
\end{align*}
\end{theorem}
\begin{proof}
Given $\epsilon < \epsilon_0$ from Theorem \ref{higherlocalregthm}, 
there exists an $R$ large so that 
\begin{align*}
\int_{D(R)} |Rm|^2 dV_g < \epsilon < \epsilon_0.
\end{align*}
Choose any $x \in X$ with $d(x,p) = r(x) > 2R$, 
then $B(x,r) \subset D(R)$. From 
Theorem \ref{higherlocalregthm}, we have 
\begin{align*}
\underset{B(x, r/2)}{sup}| \nabla^k Rm| \leq
\frac{C_k}{r^{2+k}} \left\{ \int_{B(x,r)} |Rm|^2 dV_g \right\}^{1/2}
\leq \frac{C_k \epsilon}{r^{2+k}}, 
\end{align*}
which implies
\begin{align}
\label{cond1}
| \nabla^k Rm|(x)  \leq \frac{C_k \epsilon}{r^{2+k}}.
\end{align}
\end{proof}
  Clearly, as we take $R$ larger, we may choose 
$\epsilon$ smaller, and we see that
\begin{align}
\label{imp2'}
\underset{S(r)}{sup } \ |\nabla^k Rm| = o(r^{-2-k}),
\end{align}
as $r \rightarrow \infty$, where $S(r)$ denotes the 
sphere of radius $r$ centered at $p$.
\section{Volume Growth}
\label{euclideanvolumegrowth}
This section will be devoted to proving the following
\begin{theorem}
\label{bigthm}
Let $(X,g)$ be a complete, noncompact, 4-dimensional 
Riemannian manifold with base point $p$.
Assume that there exists a constant $C_1 > 0$ so that
\begin{align}
\label{cond4}
Vol(B(q,s)) \geq C_1 s^4,
\end{align}
for any $q \in X$, and all $s \geq 0$.
Assume furthermore that as $r \rightarrow \infty$,
\begin{align}
\label{decay1}
\underset{S(r)}{sup} \ |Rm_g| &= o(r^{-2}),
\end{align}
where $S(r)$ is sphere of radius $r$ centered at $p$. 
If $b_1(X) < \infty$, then $(X,g)$ has finitely many ends, 
and there exists a constant $C_2$ (depending on $g$) so that 
\begin{align}
\label{vga}
Vol(B(p,r)) \leq C_2 r^4.
\end{align}
Furthermore, each end is ALE of order $0$. 
\end{theorem}
\begin{remark}
As stated in the introduction, this theorem 
holds in dimension $n$. For simplicity, we consider
the $4$ dimensional case - the same proof 
applies in dimension $n$ with appropriate 
modification of constants. 
\end{remark}
\begin{remark} We emphasize that our proof requires a weaker condition than 
$b_1(X) < \infty$. The condition is that there are
only finitely many disjoint ``bad'' annuli components in $X$. 
We say a component $A_0(r_1,r_2)$ of an annulus 
$A(r_1, r_2) = \{ q \in X \ | \ r_1 < d(p, q) < r_2 \}$
is {\em{bad}} if $S(r_1) \cap \overline{A_0(r_1,r_2)}$ has 
more than $1$ component, where $S(r_1)$ is the 
sphere of radius $r_1$ centered at $p$. 
If $b_1(X) < \infty$ then $X$ may contain only finitely many disjoint bad annuli
(see the Mayer-Vietoris argument in Lemma \ref{components} below).
We expect that this assumption can be removed, but at the 
moment substantial technical difficulties remain. 
\end{remark}
\subsection{Area Comparison}
 For notation, we let $\rho(x) = d(p,x)$ denote the 
distance function from $p$, and recall that $\Delta \rho (x)$
is the mean curvature of the level set of $\rho$, 
at any smooth point $x$ of $\rho$. 
We note that the distance function $\rho(x)$ is 
only Lipschitz, but $\Delta \rho(x)$ is well-defined
almost everywhere. Letting $\mathcal{C}(p)$ denote the cut 
locus of $p$, it is well-known that $\mathcal{C}(p)$ has 
measure zero, and $\rho(x)$ is smooth on $X \setminus \mathcal{C}(p)$
(see \cite[Chapter 3]{Chavel}).
\begin{lemma}
\label{radialcomp}
There exists a decreasing function $\epsilon(r) \geq 0$ with 
$\epsilon(r) \rightarrow 0$ as $ r \rightarrow \infty$,
such that 
\begin{align}
\label{Hdecay}
\Delta \rho(x) \leq  \frac{3+ \epsilon(\rho(x))}{\rho(x)},
\end{align}
at any point $x$ where the distance function $\rho(\cdot) = d(p, \cdot)$
is smooth. 
\end{lemma}
\begin{proof}
First we construct a
radial comparison metric as follows. Let $ \tilde{g}
= dr^2 + h(r)^2 d\theta^2$. 
The Jacobi equation is 
\begin{align}
\label{Ricatti}
h^{''} = - K(r) h,
\end{align} where 
$K(r)$ is the radial curvature. 
We choose a smooth radial curvature function $K(r)$ so
that 
\begin{align*}
&3K(r) \leq \mbox{min} \{ Ric_g, 0 \}\\
&3K(r) = - A r^{-2} \mbox{ for } r > r_0 \mbox{ large}.
\end{align*}
The existence of a radial metric $\tilde{g}$ satisfying 
$h(0) = 0, h'(0) = 1$, and with $K(r)$ as its 
radial curvature function is obtained easily 
(see \cite{GreeneWu}, Proposition 4.2).
We see that $r^q$ solves the Jacobi equation on $(r_0, \infty)$
where $q(q-1) = \frac{A}{3}$. We thus have two 
linearly independent solutions, and therefore the
general solution is given by  
$h(r) = c_1 r^{q_+} + c_2 r^{q_-}$ on $(r_0, \infty)$ 
for some constants $c_1, c_2$,
where $q_{\pm} = \frac{1}{2} \left( 1 \pm \sqrt{ ( 1 + 
\frac{4A}{3}} \right).$
We may arrange so that our comparison metric 
satisfies $c_2 = 0$,  and we then have for $r > r_0$,
\begin{align*}
\frac{h'}{h} = \frac{ c_1 q_+ r^{q_{+} -1}}{c_1 r^{q_{+}}} = \frac{q_{+}}{r}.
\end{align*}
Using this radial metric and the decay condition (\ref{decay1}), 
the lemma follows from the Laplacian comparison theorem 
(see \cite{GreeneWu}). 
\end{proof}
\begin{definition}
\label{hausdorff}
Let $S(r)$ denote the geodesic sphere 
of radius $r$. 
Since the distance function is Lipschitz, 
by the coarea formula, for almost every $r$, 
$S(r)$ is $\mathcal{H}^{3}$-measurable,
where $\mathcal{H}^{3}$ denotes $3-$Hausdorff 
measure, and we define 
\begin{align} 
\mathcal{H}(r) = \mathcal{H}^{3}(S(r)). 
\end{align}
For a set $K$ we let 
\begin{align}
\mathcal{H}_K(r) = \mathcal{H}^{3}(S(r) \cap K).
\end{align}
\end{definition}
Next we define the lower area of distance spheres.
\begin{definition}
\label{lowerarea}
Let $\mathbf{S}_p$ denote the unit sphere 
in $T_p(X)$, and $\mathbf{D}_p$ denote the maximal star-shaped region 
on which $\mbox{exp}_p$ is a diffeomorphism. 
Also, let $\mathbf{D}_p(r)$ be the subset 
of $\mathbf{S}_p$ of directions $\xi$ such 
that $ r \xi \in  \mathbf{D}_p$, that is, 
\begin{align}
r \cdot \mathbf{D}_p(r) = \mathbf{S}(r) \cap \mathbf{D}_p. 
\end{align}
We define the lower area of $S(r)$ by 
\begin{align}
\mathcal{A}(p,r) = \int_{\mathbf{D}_p(r)} \sqrt{ \mathbf{g}}(r;\xi) 
d \sigma,
\end{align}
where $\sqrt{ \mathbf{g}}(r;\xi)$ is given by 
$\exp^*dV_g(r;\xi) =  \sqrt{ \mathbf{g}}(r;\xi)
dr d \sigma $, and $d \sigma$ is the spherical area 
element of $\mathbf{S}_p$, induced by Lebesgue 
measure on $T_p(X)$.  
\end{definition}
The following proposition gives the relation between 
the lower area and the $3$-Hausdorff measure of distance spheres:
\begin{proposition}
\label{eqofarea}
For almost every $r$, 
\begin{align}
\mathcal{H}(r) = \mathcal{A}(r).
\end{align}
\end{proposition}
\begin{proof}
This is an easy consequence of the coarea formula, 
see \cite[Proposition 3.4]{Chavel}. 
\end{proof}
\begin{theorem}[Area Comparison]
\label{areacomp}
There exists a function $\epsilon(r)$ with 
$\epsilon(r) \rightarrow 0$ as $r \rightarrow \infty$ 
such that for almost every $r_1$ and $r_2$, 
\begin{align}
\mathcal{H}(r_2) \leq \mathcal{H}(r_1) 
\left( \frac{r_2}{r_1} \right)^{3 + \epsilon(r_1)}.
\end{align}
\end{theorem}
\begin{proof}
In spherical exponential coordinates, we have the formula \cite{Petersen}, 
\begin{align}
\partial_r \sqrt{ \mathbf{g}}(r; \theta) 
= \Delta \rho \cdot   \sqrt{ \mathbf{g}}(r; \theta).
\end{align}
Integrating, we obtain
\begin{align*}
 \sqrt{ \mathbf{g}}(r_2; \theta) = 
\sqrt{ \mathbf{g}}(r_1; \theta) \cdot
\exp{ \left( 
\int_{r_1}^{r_2} \Delta \rho (t; \theta) dt \right)}.
\end{align*}
Using Lemma \ref{radialcomp}, we obtain 
\begin{align*}
 \sqrt{ \mathbf{g}}(r_2; \theta)
 \leq 
 \sqrt{ \mathbf{g}}(r_1; \theta) \cdot
\exp{ \left( 
\int_{r_1}^{r_2} \frac{3 + \epsilon}{t} dt \right)}
\leq 
 \sqrt{ \mathbf{g}}(r_1; \theta)
\left( \frac{r_2}{r_1} \right)^{3 + \epsilon}.
\end{align*}
Recalling the notation introduced in Definition \ref{lowerarea}, 
we have that for $r_1 < r_2$, 
\begin{align}
\mathbf{D}_p (r_2) \subseteq {\mathbf{D}_p(r_1)},
\end{align}
that is, if a point $x \in \exp_p \{(r_2, \mathbf{D}_p(r_2) )\}$, then 
there exists a minimal geodesic from $p$ to $x$, and this geodesic 
will necessarily hit the distance sphere of radius $r_1$. 

 Next using Proposition \ref{eqofarea}, for almost every 
$r_1$, $r_2$, 
\begin{align*}
\mathcal{H}(r_2) = \mathcal{A}(p,r_2)  
&=\int_{ \mathbf{D}_p(r_2)} \sqrt{\mathbf{g}}(r_2, \theta)
 d \sigma\\
&\leq \int_{ \mathbf{D}_p(r_2)} 
 \sqrt{ \mathbf{g}}(r_1, \theta)
\left( \frac{r_2}{r_1} \right)^{3 + \epsilon}
 d \sigma \\
& \leq \left( \frac{r_2}{r_1} \right)^{3 + \epsilon(r_1)} 
\int_{ \mathbf{D}_p(r_1)} \sqrt{ \mathbf{g}}(r_1, 
\theta) d \sigma 
= \mathcal{H}(r_1) \left( \frac{r_2}{r_1} \right)^{3 + \epsilon(r_1)}. 
\end{align*}
\end{proof}

\subsection{Selection of Annuli}
Fix $s > 1$, and we now choose a sequence of components of annuli 
in the following manner. Start with $A_0(1,s)$ = any component 
of $A(1,s)$. The outer boundary $\partial A_0(1,s)$, which we denote 
by $S_{outer,0}$, may have several components. Choose any component, 
and call this $S_{inner,1}$. Next consider the annulus 
$A(s,s^2)$. This may have several components, but we 
define $A_1(s,s^2)$ to be the component which has inner 
boundary portion $S_{inner,1}$. We repeat this procedure:
given $A_j(s^j, s^{j+1})$, choose any component of the 
outer boundary, $S_{outer, j}$, and call this component 
$S_{inner, j+1}$. Choose $A_{j+1}(s^{j+1}, s^{j+2})$ 
to be the component of $A(s^{j+1}, s^{j+2})$ which has 
inner boundary portion $S_{inner, j+1}$.

\begin{lemma}
\label{components}
Let $A_j(s^j, s^{j+1})$ be any sequence of components 
of the metric annulus $A(s^j, s^{j+1})$. 
Then $S(s^j) \cap  \overline{A_j(s^j, s^{j+1})}$ has only 1 component,
except for finitely many $j$. 
That is, there are only finitely many $A_j(s^j, s^{j+1})$
with the initial boundary sphere having more than 1 component.
Therefore, there exists an $N$ such that for all 
$j > N$, $S(s^j) \cap  \overline{A_j(s^j, s^{j+1})}$ has only 1 component.
\end{lemma}
\begin{proof}
From the assumption, $H_1(X)$ is finitely generated, 
say there are at most $k$ generators. Assume that 
there are $k+1$ annuli with  $S(s^j) \cap  \overline{A_j(s^j, s^{j+1})}$
having more than 1 component, index these by $j_i$, 
$i = 1 \dots k+1$. 
For simplicity, assume that number of components
of the initial boundary of $\overline{A_j}$ is $2$, and that 
the number of components of the outer boundary of $\overline{A_j}$ is $1$.
Let $ U = \coprod A_j(s^j, s^{j+1})$,
and write $ X = U \cup V$, where $V$ is a open 
set which intersects $A_j$ in an $\epsilon$-neighborhood 
of each boundary component of $A_j$.
Notice that $U$ has $k+1$ components, and by choosing 
$\epsilon$ sufficiently small, $U \cap V$ 
has $ 3(k+1)$ components. The important observation 
is that $V$ has at most $k+2$ components.
This is because all of the initial spheres 
are connected to the base point by a geodesic, 
so this gives $1$ connected component of $V$,
and there are at most $k+1$ more components 
of $V$ which touch each outer boundary sphere. 
We then consider the following portion of the 
Mayer-Vietoris sequence in homology:
\begin{align*}
&H_1(X) \rightarrow H_0(U \cap V) \rightarrow
H_0(u) \oplus H_0(V) \rightarrow H_0(X) \rightarrow 0,
\end{align*}
By the observation above on the number of components
of $U$, $V$, and $U \cap V$, this sequence is
\begin{align*}
H_1(X) \rightarrow \mathbb{Z}^{3k} \rightarrow
\mathbb{Z}^k \oplus \mathbb{Z} \oplus\mathbb{Z}^{\leq k+1} 
\rightarrow \mathbb{Z} \rightarrow 0. 
\end{align*}
It is then easy to see this forces $H_1(X)$ to contain 
a $\mathbb{Z}^{k+1}$, which contradicts the fact that 
$H_1(X)$ has at most $k$ generators. 
For simplicity we have restricted to this 
simple case, but a similar argument shows that 
if we have $k+1$ disjoint annuli with the 
initial boundary of each having more than 
1 component, then this forces $H_1(X)$ to 
contain at least a $\mathbb{Z}^{k+1}$.
\end{proof}
\begin{proposition}
\label{annchoose}
There exists a subsequence $\{j \} \subset \{ i \}$  satisfying 
\begin{align}
\label{annulicondo}
\mathcal{H}^3( S_{inner, j+1} ) \geq  (1 - \eta_j) 
\mathcal{H}^3( S_{inner, j} )s^3,
\end{align}
where $\eta_j \rightarrow 0$ as $j \rightarrow \infty$.
\end{proposition}
\begin{proof} Consider the sequence 
of components of annuli $\{A_i(s^i, s^{i+1}) \}$ for 
$ i = 0 \dots \infty$. 
If there does not exists such a subsequence satisfying
(\ref{annulicondo}), then for all $j \geq N$, we have 
\begin{align}
\mathcal{H}^3( S_{inner, j+1}) \leq(1 - \eta') 
\mathcal{H}^3( S_{inner, j}) s^3,
\end{align}
and $\eta' > 0 $. 
Given any $r \geq s^N$, choose $j$ so that 
$$s^j \leq r \leq s^{j+1}.$$

Theorem \ref{areacomp} implies that 
there exists a function $\epsilon(r)$ with 
$\epsilon(r) \rightarrow 0$ as $r \rightarrow \infty$ 
such that 
\begin{align}
\mathcal{H}^3(S(r) \cap A_j) \leq \mathcal{H}^3(S_{inner,j}) s^{3+ \epsilon}
\end{align}
Note we are using Lemma \ref{components}, since the 
initial boundary has only $1$ component for all 
$j$ large, we can apply the area comparison. 
Note also that Theorem \ref{areacomp} holds for almost every 
$r_1$, $r_2$, but in the above and what follows, 
we are free to change the sequence chosen 
by arbitrarily small changes, and we will do 
this automatically. 

Since $s$ is bounded, and $\epsilon \rightarrow 0$ 
as $r \rightarrow \infty$, we have
\begin{align}
\mathcal{H}^3(S(r) \cap A_j ) \leq \mathcal{H}^3(S_{inner,j}) 
s^{3+ \epsilon} \leq \mathcal{H}^3(S_{inner,j}) 
s^3 (1 + \epsilon).
\end{align}
Therefore 
\begin{align*}
\mathcal{H}^3(S(r) \cap A_j) &\leq ( 1 + \epsilon) s^3 
(1 - \eta')\mathcal{H}^3(S_{inner,j-1}) s^3\\
&\leq  ( 1 + \epsilon) s^3  (1 - \eta')^{j-N}
\mathcal{H}^3(S_{inner, N}) s^{3(j-N)}\\
&=  ( 1 + \epsilon)(1 - \eta')^{j-N} 
s^{3-3N}  \mathcal{H}^3(S_{inner,N}) s^{3j}\\
& \leq C (1-\eta')^{j-N}s^{3j} \leq C (1 - \eta')^{j-N} r^3.
\end{align*}
From the coarea formula, we have for $j >>N$,  
\begin{align*}
Vol( A_j(s^j, s^{j+1})) = \int_{s^j}^{s^{j+1}} \mathcal{H}( S(t) \cap
A_j) dt \leq C(1- \eta')^{j-N}( s^{4(j+1)} - s^{4j}).
\end{align*}
But the condition (\ref{cond4}) implies that 
there exists $C' > 0$ so that  
\begin{align*}
Vol( A_j(s^j, s^{j+1})) \geq C' s^{4(j+1)},
\end{align*}
a contradiction. 
\end{proof}
 In the following we will take our subsequence $\{ j \}$
to be maximal, that is, so that 
\begin{align}
\mathcal{H}^3( S_{inner, i+1} ) <  \mathcal{H}^3( S_{inner, i} )s^3,
\end{align}
for $i$ not in our subsequence $\{j \} \subset \{i\}$.
This is possible since if 
\begin{align*}
\mathcal{H}^3( S_{inner, i+1} ) \geq  \mathcal{H}^3( S_{inner, i} )s^3,
\end{align*}
then we may obviously include $i$ in the subsequence $\{j\}$. 
In the following we will reserve the index 
$j$ for this subsequence, while the index $i$ will 
index all annuli. 
\subsection{Laplacian of the distance function}
Some of the material in this section in well-known to 
experts, but we for completeness, we include the proof. 

\begin{lemma}
\label{intlap1} For almost every $r_1, r_2$, 
\begin{align}
\int_{A(r_1, r_2)} (- \Delta \rho) dV_g \leq \mathcal{H}(r_1) - 
\mathcal{H}(r_2)
\end{align}
\end{lemma}
\begin{remark} 
In the above, we are viewing $\Delta \rho$ as a measurable 
function, defined almost everywhere. 
\end{remark}
\begin{proof}
We begin with an approximation process as in 
\cite[Proposition 1.1]{SchoenYau2}.
As before, we let $\mathbf{D}_p$ denote the maximal star 
shaped domain inside of the $T_p(X)$ such that 
$exp : \mathbf{D}_p \mapsto exp_p(E_p) =D_p$ is 
a diffeomorphism, and $\mathcal{C}(p) = \partial D_p$,
where $\mathcal{C}(p)$ is the cut locus. We note that $\mathcal{C}(p)$ 
has measure zero and $X = D_p \cup \mathcal{C}(p)$.
The distance function $\rho(x) = d(p,x)$ is 
Lipschitz, smooth on $X \setminus \mathcal{C}(p)$ and satisfies
$|\nabla \rho |^2 = 1$ on $X \setminus \{ \mathcal{C}(p) \cup p \}$.
We claim that for any $\phi \in C^{\infty}_{c}(X)$,
$\phi \geq 0$, with support of $\phi$ not containing $p$, we have 
\begin{align}
\label{intparts}
- \int_X \phi \Delta \rho \leq \int_X \nabla \phi \cdot 
\nabla r.
\end{align}
To prove this, let $D_p = exp_p(E_p)$, and since the 
cut locus has measure zero, we have 
\begin{align*}
\int_X \phi \Delta \rho = \int_{D_p} \phi \Delta \rho.
\end{align*}
Since $\mathbf{D}_p$ is star shaped, we may construct a family of 
smooth star-shaped domains $D_{\epsilon} \subset 
D_p$ with $\lim_{\epsilon \rightarrow 0} D_{\epsilon}
= D_p$. Since $D_{\epsilon}$ is star-shaped, 
we have
\begin{align*}
\frac{\partial r}{\partial \nu_{\epsilon}} > 0 \mbox{ on }
\partial D_{\epsilon}. 
\end{align*}
On $D_{\epsilon}$, $r$ is smooth, so applying 
Green's first identity, we have 
\begin{align*}
\int_{D_{\epsilon}} \phi \Delta \rho
=  - \int_{D_{\epsilon}} \nabla \phi \cdot 
\nabla \rho
+ \int_{ \partial D_{\epsilon}} \phi 
\frac{\partial r}{\partial \nu_{\epsilon}}
\geq - \int_{D_{\epsilon}} \nabla \phi \cdot 
\nabla \rho.
\end{align*}
Since $|\nabla \rho| = 1$ almost everywhere, by dominated 
convergence we have
\begin{align*}
\limsup_{\epsilon \rightarrow 0}
\int_{D_{\epsilon}} \phi \Delta \rho
\geq  - \int_{D_p} \nabla \phi \cdot 
\nabla \rho =  - \int_{X} \nabla \phi \cdot 
\nabla \rho
\end{align*}
Next, from Lemma \ref{radialcomp}, $\Delta \rho \leq C$ 
almost everywhere on the complement of any 
compact subset containing $p$, so from the 
assumption on the support of $\phi$, we have 
$C - \Delta \rho \geq 0$. By Fatou's Lemma, we conclude that 
\begin{align}
\label{fatou}
 \int_{D_p} \phi(C - \Delta \rho ) \leq \liminf_{\epsilon \rightarrow 0}
\int_{D_{\epsilon}} \phi (C -\Delta \rho).
\end{align}
This yields 
\begin{align*}
 \int_{X} \phi \Delta \rho
= \int_{D_p} \phi \Delta \rho  \geq \limsup_{\epsilon \rightarrow 0}
\int_{D_{\epsilon}} \phi \Delta \rho \geq - \int_{X} \nabla \phi \cdot 
\nabla \rho.
 \end{align*}
which proves (\ref{intparts}). Note that by an approximation 
argument, (\ref{intparts}) holds for $\phi \in C^{0,1}_{c}(X)$, 
that is Lipschitz functions with compact support. 

 Next, define Lipschitz cutoff functions 
$a_{\epsilon} \in C^{0,1}(r_1,r_2)$ such that 
\begin{align}
a_{\epsilon} 
= \left\{ 
\begin{array}{ll} 
1 & \mbox{ if } r_1 + \epsilon \leq t \leq r_2 - \epsilon \\
0 & \mbox{ if } t \leq r_1 \mbox{ or } t \geq r_2
\end{array} \right.
\end{align}
and such that $|\nabla a_{\epsilon}| = \epsilon^{-1}$ for 
$r_1 \leq t \leq r_1 + \epsilon$ and 
$r_2 - \epsilon \leq t \leq r_2$. 
Let $\phi_{\epsilon} (x)= a_{\epsilon}(\rho(x))$, 
which is also a Lipschitz function, therefore we may apply 
(\ref{intparts}) to obtain 
\begin{align*}
-\int_{A(r_1,r_2)} \phi_{\epsilon} \Delta \rho 
&= -\int_{X} \phi_{\epsilon} \Delta \rho 
\leq \int_X \nabla \phi_{\epsilon} \cdot \nabla \rho\\
& = \int_{A(r_1, r_1 + \epsilon)} 
\frac{1}{\epsilon}  - \int_{A(r_2 - \epsilon, r_2)} 
\frac{1}{\epsilon}\\
& = \frac{ Vol(A(r_1, r_1 + \epsilon))}{\epsilon}
- \frac{ Vol (A(r_2 - \epsilon, r_2))}{\epsilon}.
\end{align*}
Since $\Delta \rho$ is bounded from above away from $p$,
an application of Fatou's Lemma as in (\ref{fatou}) 
above, yields the inequality 
\begin{align}
- \int_{A(r_1,r_2)} \Delta \rho  \leq \liminf_{\epsilon \rightarrow 0}
- \int \phi_{\epsilon} \Delta \rho.
\end{align}
From the coarea formula 
\begin{align}
 Vol(A(0,t)) = \int_0^t \mathcal{H}(s) ds,
\end{align} 
therefore $Vol(A(0,t))$ is differentiable almost everywhere 
with derivative $\mathcal{H}(t)$. Since 
$$Vol(A(a,t)) = Vol(A(0,t)) - Vol(A(0,a)),$$ we have 
that $Vol(A(a,t))$ is differentiable at $t=a$ for almost 
every $a$ with derivative $\mathcal{H}(a)$. Similarly, 
$Vol(A(t,a))$ is differentiable at $t=a$ for almost 
every $a$ with derivative $\mathcal{H}(a)$.
Taking limits in the above, we have 
\begin{align}
\label{gogog}
- \int_{A(r_1,r_2)} \Delta \rho  \leq
\mathcal{H}(r_1) - \mathcal{H}(r_2),
\end{align}
for almost every $r_1$ and $r_2$. 
\end{proof}
\begin{corollary} 
\label{deltarl1}
\begin{align}
\Delta \rho \in L^1_{loc}.
\end{align}
\end{corollary}
\begin{proof}
As $x \rightarrow p$, $\Delta \rho \sim \frac{3}{r}$, so the 
result is true for any small enough ball around $p$.
On any annulus $A(r_1,r_2)$, $r_1 >0$, by  Lemma \ref{radialcomp},
$\Delta \rho \leq C$ almost everywhere. Therefore
for almost every $r_1, r_2$  
\begin{align*}
\int_{A(r_1,r_2)} |\Delta \rho|
&= \int_{A(r_1,r_2)} | \Delta \rho -C + C| 
\leq  C Vol({A(r_1,r_2)}) + \int_{A(r_1,r_2)} |\Delta \rho - C|\\
&= C Vol({A(r_1,r_2)}) + \int_{A(r_1,r_2)} (C - \Delta \rho) \\
& \leq  2 C Vol({A(r_1,r_2)}) + \mathcal{H}(r_1) - \mathcal{H}(r_2).
\end{align*}
\end{proof} 
Next we want the explicit error term in Lemma \ref{intlap1}. 
We recall the following from \cite{ItohTanaka}, \cite{Ozols}, 
\cite{Mantegazza}. The cutlocus is $(n-1)$-rectifiable and can be 
decomposed as the disjoint union 
$\mathcal{C} = \mathcal{C}_1 \cup \mathcal{C}_2$, 
where $\mathcal{C}_2$ has finite $\mathcal{H}^{n-2}$-measure, 
and $\mathcal{C}_1$ is a union of smooth hypersurfaces. 
The set $\mathcal{C}_1$ is characterized 
by the following: If $x \in \mathcal{C}_1$ then 
(i) $x$ is not conjugate to $p$, and 
(ii) there are exactly $2$ geodesics joining $p$ to $x$. 
The set $\mathcal{C}_2$ is characterized by: 
if $x \in \mathcal{C}_1$ then $x$ is either a 
conjugate point, or there are more than $2$ minimal geodesics 
from $p$ to $x$. 
\begin{lemma}
\label{intlap3} 
For almost every $r_1, r_2$, 
\begin{align}
\int_{A(r_1, r_2)} (- \Delta \rho) dV_g = \mathcal{H}(r_1) - 
\mathcal{H}(r_2) 
- \int_{\mathcal{C}}
   \Big(  \langle (\nabla \rho)^{+}, \partial_{\mathcal{C}}^{+} 
\rangle + \langle (\nabla \rho)^{-}, \partial_{\mathcal{C}}^{-} 
\rangle \Big) d \mathcal{H}^3,
\end{align}
Since $\mathcal{C}_2$ has $\mathcal{H}^{3}$-measure zero, 
the error term integration is only over 
$\mathcal{C}_1$, and there are 2 terms corresponding to the 
2 radial directions, with $\partial_{\mathcal{C}}^{\pm}$ 
denoting the corresponding outward normals to $\mathcal{C}$.
\end{lemma}
\begin{proof}
 We first claim that 
\begin{align}
\label{intparts2}
- \int_X \phi \Delta \rho = \int_X \nabla \phi \cdot 
\nabla \rho -  \int_{\mathcal{C} \cap A(r_1,r_2)}
   \phi \Big(  \langle (\nabla \rho)^{+}, \partial_{\mathcal{C}}^{+} 
\rangle + \langle (\nabla \rho)^{-}, \partial_{\mathcal{C}}^{-} 
\rangle \Big) d \mathcal{H}^3,  
\end{align}
To prove this, we follow the proof of Lemma \ref{intlap1}, 
applying Green's identity, we have
\begin{align*}
\int_{D_{\epsilon}} \phi \Delta \rho
=  - \int_{D_{\epsilon}} \nabla \phi \cdot 
\nabla \rho
+ \int_{ \partial D_{\epsilon}} \phi 
\frac{\partial r}{\partial \nu_{\epsilon}}.
\end{align*}

From Lemma \ref{deltarl1}, $\Delta \rho \in L^1_{loc}$ so 
we may apply dominated convergence to show that
\begin{align*}
- \int_X \phi \Delta \rho &= \lim_{\epsilon \rightarrow 0}
\int_{D_{\epsilon}} \phi \Delta \rho
= - \int_{D_p} \phi \Delta \rho \\
& = \lim_{\epsilon \rightarrow 0} \Big( 
- \int_{D_{\epsilon}} \nabla \phi \cdot \nabla \rho
+ \int_{ \partial D_{\epsilon}} \phi 
\frac{\partial r}{\partial \nu_{\epsilon}} \Big)\\
& = - \int_{X} \nabla \phi \cdot \nabla \rho
-  \lim_{\epsilon \rightarrow 0} \Big( \int_{ \partial D_{\epsilon}} \phi 
\frac{\partial r}{\partial \nu_{\epsilon}} \Big).
\end{align*}
From the remarks above, away from a set of finite $\mathcal{H}^2$-measure,
the cut locus is a union of smooth hypersurfaces. It is then clear that 
\begin{align*}
 \lim_{\epsilon \rightarrow 0}\frac{\partial r}{\partial \nu_{\epsilon}}
= \langle \partial_{\mathcal{C}}, \nabla \rho \rangle,
\end{align*}
except on a set of $\mathcal{H}^3$-measure zero. So we have 
\begin{align}
\label{errtermdef}
  \lim_{\epsilon \rightarrow 0} \Big( \int_{ \partial D_{\epsilon}} \phi 
\frac{\partial r}{\partial \nu_{\epsilon}} \Big)
=-  \int_{\mathcal{C}}
   \phi \Big(  \langle (\nabla \rho)^{+}, \partial_{\mathcal{C}}^{+} 
\rangle + \langle (\nabla \rho)^{-}, \partial_{\mathcal{C}}^{-} 
\rangle \Big) d \mathcal{H}^3,
\end{align}
which proves (\ref{intparts2}). Finally we may imitate 
the argument in Lemma \ref{intlap1}, using the sequence 
of cutoff functions $a_{\epsilon}$, and the error 
term is
\begin{align*}
&- \lim_{\epsilon \rightarrow 0} \int_{\mathcal{C}}
   a_{\epsilon} \Big(  \langle (\nabla \rho)^{+}, \partial_{\mathcal{C}}^{+} 
\rangle + \langle (\nabla \rho)^{-}, \partial_{\mathcal{C}}^{-} 
\rangle \Big) d \mathcal{H}^3\\
& = - \int_{\mathcal{C} \cap A(r_1,r_2)} 
\Big(  \langle (\nabla \rho)^{+}, \partial_{\mathcal{C}}^{+} 
\rangle + \langle (\nabla \rho)^{-}, \partial_{\mathcal{C}}^{-} 
\rangle \Big) d \mathcal{H}^3,
\end{align*}
the limit being justified since by the choice 
of $r_1$ and $r_2$, the cut locus hits 
$S(r_1)$ and $S(r_2)$ in a set of 
$\mathcal{H}^3$-measure zero. 
\end{proof}
\begin{lemma}
\label{intlap2} 
For almost every $r_1, r_2$, 
\begin{align}
\begin{split}
 - \int_{A(r_1, r_2)} \Delta (\rho^2) dV_g  & =
 2 \big( r_1 \mathcal{H}(r_1) - r_2 \mathcal{H}(r_2) \big)\\
& - 2 \int_{\mathcal{C} \cap A(r_1,r_2)} 
\Big(  \langle \rho (\nabla \rho)^{+}, \partial_{\mathcal{C}}^{+} 
\rangle + \langle \rho (\nabla \rho)^{-}, \partial_{\mathcal{C}}^{-} 
\rangle \Big) d \mathcal{H}^3,
\end{split}
\end{align}
\end{lemma}
\begin{proof}
The formula (\ref{intparts2}) is replaced with  
\begin{align}
\label{intparts3}
- \int_X \phi \Delta \rho^2 = 2 \int_X \nabla \phi \cdot 
\rho \nabla \rho -  2 \int_{\mathcal{C} \cap A(r_1,r_2)}
   \phi \Big(  \langle  \rho (\nabla \rho)^{+}, \partial_{\mathcal{C}}^{+} 
\rangle + \langle \rho (\nabla \rho )^{-}, \partial_{\mathcal{C}}^{-} 
\rangle \Big) d \mathcal{H}^3,  
\end{align}
and the proof proceeds as in Lemma \ref{intlap3}.
\end{proof}
\begin{remark}
Our proof does not require such an explicit expression for
the error term. An examination of the proof below 
shows we just require 2 properties:
(i) that the error term is additive under finite
disjoint unions, and that (ii) the error 
term with $\phi$ can be estimated by
$ \Vert \phi \Vert_{L^{\infty}}$ times
the error term with $\phi \equiv 1$. Both of these properties 
are clear from the definition of the error term as the limit 
in (\ref{errtermdef}). 
\end{remark}
\subsection{$L^1$ convergence}

For any $s > 0$, let $(\tilde{A}^r(1, s), \tilde{g}) = 
(A(r, sr), g/r^2)$ denote the annulus $\{ x \in X : 
r < d(p, x) < sr \}$,
with the rescaled metric $g/r^2$. From (\ref{decay1})
we have 
\begin{align}
\label{decay11}
\underset{\tilde{A}^r}{\mbox{sup}}\tilde{|Rm|} \leq \epsilon_1(r).
\end{align} 
where $\epsilon_1(r) \rightarrow 0$ as $r \rightarrow \infty$.

Let $(\tilde{A}_i(1, s), \tilde{g}) = 
(A_i(s^i, s^{i+1}), g/s^{2i})$ denote our previously chosen 
sequence of annular components, but with the rescaled metric $g/s^{2i}$,
and let $\tilde{\rho}_i = s^{-i} \cdot \rho$ denote the 
rescaled distance function.  

From (\ref{decay1}) we have 
\begin{align}
\label{decay12}
\underset{\tilde{A}_i(1,s)}{\mbox{sup}}\tilde{|Rm|} \leq \epsilon_2(i), 
\end{align} 
where $\epsilon_2(i) \rightarrow 0$ as $i \rightarrow
\infty$. 
\begin{proposition}
\label{l1con}
For the subsequence $\{ j \} \subset \{i\}$,  
\begin{align}
\frac{1}{ Vol( \tilde{A}_j(1,s))}\int_{\tilde{A}_j(1,s)} \Big| 
\Delta (\tilde{\rho}_j)^2 - 8 \Big|
\tilde{dV_j} \rightarrow 0,
\end{align}
and 
\begin{align} 
\label{errorterm}
\frac{1}{Vol( \tilde{A}_j(1,s))}\int_{\mathcal{C} \cap \tilde{A}(1,s)} 
\Big(  \langle \tilde{\rho}_j (\nabla \tilde{\rho}_j)^{+}, \partial_{\mathcal{C}}^{+} 
\rangle + \langle \tilde{\rho}_j (\nabla \tilde{\rho}_j)^{-}, \partial_{\mathcal{C}}^{-} 
\rangle \Big) d \mathcal{H}^3 \rightarrow 0, 
\end{align}
as  $j \rightarrow \infty$.
\end{proposition}
\begin{proof}
For notation, we set 
\begin{align*}
E_j = \int_{\mathcal{C} \cap \tilde{A}(1,s)} 
\Big(  \langle \tilde{\rho}_j (\nabla \tilde{\rho}_j)^{+}, \partial_{\mathcal{C}}^{+} 
\rangle + \langle \tilde{\rho}_j (\nabla \tilde{\rho}_j)^{-}, \partial_{\mathcal{C}}^{-} 
\rangle \Big) d \mathcal{H}^3
\end{align*}
If the proposition is not true, then for some subsequence 
(which for simplicity we continue to index by $j$), 
\begin{align*}
\frac{1}{ Vol( \tilde{A}_j(1,s))}
\Big( \int_{ \tilde{A}_j(1,s)} \Big| \Delta (\tilde{\rho}_j)^2 - 8
\Big| \tilde{dV_j} + 2 E_j \Big) \geq C > 0,
\end{align*}
which in the original regions is 
\begin{align*}
\frac{1}{ Vol(A_j)} \Big( \int_{A_j} \Big| \Delta (\rho^2) - 8 \Big| dV_g 
+ 2 E_j \Big) \geq C > 0.
\end{align*}
Note from (\ref{radialcomp}), we have 
\begin{align}
\label{radialcomp2} \Delta (\rho^2) = 2 ( \rho \Delta \rho + |\nabla \rho|^2)
\leq 2 ( \rho \frac{ 3 + \epsilon_j}{\rho} + 1) 
= 8 + \epsilon_j',
\end{align}
at any smooth point of $\rho$. 
Letting $\mathcal{C}$ denote the cut locus,  recall 
that $\mathcal{C}$ is a set of measure zero,
and the distance function is smooth on 
$X \setminus \mathcal{C}$, so in all the following
integrals we may disregard the cut locus. 
 Using (\ref{radialcomp2}), we have 
\begin{align*}
\frac{1}{ Vol(A_j)} \int_{A_j} & \Big| \Delta (\rho^2) 
- 8 \Big| dV_g
=\frac{1}{ Vol(A_j)} \int_{A_j} \Big| \Delta (\rho^2) - (8 + \epsilon_j)
+  \epsilon_j  \Big| dV_g\\
&\leq \frac{1}{ Vol(A_j)} \Big( \int_{A_j} \Big| \Delta (\rho^2) 
- (8 + \epsilon_j) \Big|dV_g \Big) + \epsilon_j\\
& =  \frac{1}{ Vol(A_j)} \Big( \int_{A_j} ( 8 + \epsilon_j -
\Delta (\rho^2)) dV_g \Big) + \epsilon_j\\
& =  - \frac{1}{ Vol(A_j)} \Big( \int_{A_j} 
\Delta (\rho^2) dV_g \Big) + 8 + 2\epsilon_j.
\end{align*}
where $\epsilon_j \rightarrow 0$ as $j \rightarrow \infty$. 
Using Lemma \ref{intlap2}, we obtain
\begin{align*}
C \leq \frac{2}{ Vol(A_j)} \Big( s^j \cdot \mathcal{H}_{A_j}(s^j) - 
s^{j+1} \cdot \mathcal{H}_{A_j}(s^{j+1})\Big) + 8 + 2\epsilon_j.
\end{align*}
In the rescaled form, this is
\begin{align*}
 \frac{2}{ Vol(\tilde{A}_j)} 
\Big(  s \cdot \mathcal{H}_{\tilde{A}_j}(s) -\mathcal{H}_{\tilde{A}_j}(1) \Big)
\leq - C  + 8 + 2 \epsilon_j.
\end{align*}
Since $\tilde{A}_j$ are chosen to satisfy $\mathcal{H}_{\tilde{A}_j}(s) 
\geq (1 - \eta_j)\mathcal{H}_{\tilde{A}_j}(1)s^3$,
we have 
\begin{align}
\label{zagaa}
 \frac{2}{Vol(\tilde{A}_j)} \cdot  \mathcal{H}_{\tilde{A}_j}(1)
\cdot \Big( (1 - \eta_j)s^4 - 1 \Big) \leq - C  + 8 + 2 \epsilon_j.
\end{align}

From Theorem \ref{areacomp}, we have the estimate 
\begin{align*}
\mathcal{H}_{\tilde{A}_j}(t) 
\leq \mathcal{H}_{\tilde{A}_j}(1){t}^{3 + \epsilon_j}.
\end{align*}
Since $s$ is bounded, we write this as
\begin{align*}
\mathcal{H}_{\tilde{A}_j}(t) 
\leq\mathcal{H}_{\tilde{A}_j}(1){t}^{3}(1 + \epsilon_j'),
\end{align*}
for $1 \leq t \leq s$, where $\epsilon_j' \rightarrow 0$ as $i \rightarrow
\infty$. 

From the coarea formula
\begin{align}
\label{zaga}
Vol( \tilde{A}_j) &= \int_1^s \mathcal{H}_{\tilde{A}_j}(t) dt
\leq \int_1^s \mathcal{H}_{\tilde{A}_j}(1){t}^{3}(1 + \epsilon_j')dt
= \frac{1 + \epsilon_j'}{4} 
\mathcal{H}_{\tilde{A}_j}(1)(s^4 -1).
\end{align}
Substituting (\ref{zaga}) into (\ref{zagaa}), we obtain 
\begin{align}
\label{zagab}
 \frac{8}{ (1 + \epsilon_j') } 
\cdot \frac{(1 - \eta_j)s^4 - 1}{s^4 - 1} - 8 - 2 \epsilon_j
\leq - C < 0.
\end{align}
The left hand side approaches $0$ as 
$j \rightarrow \infty$. Since the right hand side is a strictly 
negative constant, this is a contradiction.
\end{proof}
\subsection{Completion of proof}
We next claim that there exists a constant 
$C$ so that $Vol( \tilde{A}_j) < C$.
To prove this, assume by contradiction that  
$Vol( \tilde{A}_j)$ is not bounded. Then 
there exists a subsequence $\{j\} \subset \{i\}$
such that  $Vol( \tilde{A}_j) \rightarrow \infty$
as $ j \rightarrow \infty$. We next find 
connected subsets with large, but bounded volume.
Intuitively, it is obvious that if the 
volume of an annulus is very large, we may cut 
the annulus into several subsets with large, but bounded volume.
But we need the cutting to have certain
non-collapsing properties, so the cutting
is rather delicate. 

\begin{proposition}
\label{kchoose}
We can find disjoint sets $ \tilde{K}_{j,i} \subset  \tilde{A}_j,
 i = 1 \dots N_j$ such that 

(i) $\tilde{K}_{j,i}$ is connected. 

(ii)  $\tilde{K}_{j,i}$ is the union of balls of radius $1/4$ centered 
at points in $\tilde{A}_j( 3/2, s - 1/2 )$.

(iii) $\frac{\pi^2}{2}(s^4-1) + 1 < Vol( \tilde{K}_{j,i} ) < C_2 $, 
with $C_2$ a bounded constant. 

(iv)  $Vol ( \tilde{A}_j \setminus \{ \cup_{l=1}^{N_j} \tilde{K}_{j,i} \})
\leq C  Vol( \cup_{l=1}^{N_j} \tilde{K}_{j,l})$
with $C$ a uniformly bounded constant. 

\end{proposition}
\begin{proof}

 We first need to shrink the annuli, 
so consider $\tilde{A}_j( 3/2, s - 1/2 )$.
This set may have several components, this 
will be dealt with below. 
The sets we will choose below will be 
unions of balls of radius $1/4$ centered 
at points in $\tilde{A}_j( 3/2, s - 1/2 )$. 
The following lemma says that the volume of the annulus 
cannot be concentrated near the beginning portion 
or end portion of the annulus. 
\begin{lemma}
\begin{align}
Vol ( \tilde{A}_j(1,s) \setminus \tilde{A}_j( 3/2, s - 1/2 ))
\leq C \cdot Vol ( \tilde{A}_j( 3/2, s - 1/2 )).
\end{align}
\end{lemma}
\newcommand{\mh}{\mathcal{H}^3}
\newcommand{\A}{\tilde{A}_j}
\begin{proof}
From the coarea formula and area comparison
(Theorem \ref{areacomp}), we have
\begin{align}
\begin{split}
\label{porl-4}
Vol ( \A (1, 3/2)) = \int_{1}^{3/2} \mh_{\tilde{A}_j}(t) dt
& \leq  \int_{1}^{3/2} \mh_{\tilde{A}_j}(1)  t^{3 + \epsilon}dt
 \leq C \cdot  \mh (\tilde{A}_j \cap S(1)) .
\end{split}
\end{align}
Using the coarea formula, we may estimate 
\begin{align} 
\begin{split}
\label{porl-3}
Vol ( \A (3/2, s-1/2)) &=
\int_{3/2}^{s-1/2} \mh_{\tilde{A}_j}(t) dt
= \int_{3/2}^{s-1/2} \frac{\mh_{\tilde{A}_j}(t)}{t^3} \cdot t^3 dt\\
& \geq  \frac{\mh_{\tilde{A}_j}(t_0)}{t_0^3}
 \int_{3/2}^{s-1/2} t^3 dt\\
& \geq C \cdot  \mh ( \tilde{A}_j \cap S(t_0)),
\end{split}
\end{align}
where $t_0$ is chosen so that 
\begin{align*}
\underset{3/2 \leq t \leq s-1/2}{\mbox{ess. inf.}}
\frac{\mathcal{H}_{\tilde{A}_j}(t)}{t^3}
= \frac{ \mathcal{H}_{\tilde{A}_j}(t_0)}{t_0^3}.
\end{align*} 
That is, $t_0$ is the minimal area sphere (when rescaled to 
unit size) for $3/2 \leq t \leq s-1/2$.
If the minimum value is not actually achieved
(i.e, in case $ \mathcal{H}_{\tilde{A}_j}(t)t^{-3}$ has 
a discontinuity), then we approximate by a sequence approaching 
the minimum, we omit the details. 

From area comparison
and our choice of annulus, we have
\begin{align}
\begin{split}
\label{porl-2}
 (1 - \eta_j) \mh ( \tilde{A}_j \cap S(1))s^3 
&\leq \mh ( S_{inner, j+1}) \\
& \leq \mh ( \tilde{A}_j \cap S(s)) 
\leq C \cdot \mh ( \tilde{A}_j \cap S( t_0)).
\end{split}
\end{align}
Inequalities (\ref{porl-4}), (\ref{porl-3}), and (\ref{porl-2}) yield
\begin{align}
\label{porl0}
 Vol ( \A(1, 3/2)) \leq C Vol ( \A ( 3/2, s - 1/2)).
\end{align}
Next, from coarea and area comparison,
\begin{align}
\label{porl1}
Vol ( \A (s-1/2, s)) \leq C \mh ( \tilde{A}_j \cap S( s-1/2)).
\end{align}
Using area comparison,
\begin{align}
\label{porl3}
\mh ( \tilde{A}_j \cap S (s-1/2)) \leq C \mh( \tilde{A}_j \cap S(t_0)). 
\end{align}
Using the inequalities (\ref{porl-3}),(\ref{porl1}),  and (\ref{porl3}),
we obtain
\begin{align}
\begin{split}
\label{porl4}
Vol ( \A (s-1/2, s)) & \leq C \mh ( \tilde{A}_j \cap S( s-1/2))\\
& \leq  C \mh( \tilde{A}_j \cap S(t_0)) \leq C Vol ( \A (3/2, s-1/2)).
\end{split}
\end{align}
Inequalities (\ref{porl0}) and (\ref{porl4}) together yield
\begin{align}
Vol (  \A(1, 3/2)) \cup  \A (s-1/2, s))
\leq C \cdot Vol ( \A (3/2, s-1/2)).
\end{align}
\end{proof}

  A technical point, we started with a connected 
component $\tilde{A}_j( 1, s )$, but after shrinking, 
$\tilde{A}_j( 3/2, s - 1/2 )$ may
have several components. The next lemma 
says that we only need consider 
the component which goes ``in the 
direction'' of $S_{inner, j+1}$. 
\begin{lemma}
$\tilde{A}_j( 3/2, s - 1/2 )$ has 
exactly 1 component $\tilde{A}_j^{\circ}( 3/2, s - 1/2 )$
in the direction of $S_{inner, j+1}$, in the sense that 
any radial geodesic which hits this component, and 
lasts until $r = s$, must hit the outer 
spherical portion $S_{inner, j+1}$. For this component we have the 
estimate 
\begin{align}
Vol ( \tilde{A}_j( 3/2, s - 1/2 )
\setminus  \tilde{A}_j^{\circ}( 3/2, s - 1/2 ))
\leq C \cdot Vol ( \tilde{A}_j^{\circ}( 3/2, s - 1/2 )).
\end{align}
\end{lemma}

\begin{proof}
 From Lemma \ref{components} our selection of annuli have the property 
that $S_{inner, j}$ has one component, but $S_{outer, j}$
may have several components. 
The component in the direction 
of the next annulus we have labeled $S_{inner, j+1}$. 
When we shrink the annulus, we may have 
many components. If $2$ of these components 
have radial geodesics which extend to $S_{inner, j+1}$,
then in an argument similar to the proof of Lemma \ref{components}, 
this would yield another generator of $H_1(X)$. But since there 
are only finitely many generators of $H_1(X)$, this 
cannot happen for $R$ sufficiently large. 

 For the volume estimate, we look more closely at the 
proof of the Area Comparison Theorem \ref{areacomp}. 
In the last step, we can make an improvement, namely 
instead of integrating over $\mathbf{D}_p(r_1)$, we need only 
integrate over the directions in $\mathbf{D}_p(r_1)$ whose 
corresponding geodesics do not hit the cut 
locus before reaching $r_2$, which is exactly $\mathbf{D}_p (r_2)$ 
from the inclusion  $\mathbf{D}_p(r_2) \subset \mathbf{D}_p (r_1)$.
This yields the 
improved estimate
\begin{align}
\mathcal{H}(r_2) \leq \mathcal{H}^3 \Big\{ \mbox{exp} \Big( r_1 
\cdot \mathbf{D}_p (r_2) \Big) \Big\} \left( \frac{r_2}{r_1}
\right)^{3 + \epsilon(r_1)}. 
\end{align}
Our choice of annulus satisfies 
$\mh (S_{inner, j+1} ) \geq (1 - \eta_j)\mh(S_{inner, j})s^3$, so 
together with the above, after rescaling we have 
\begin{align}
\label{gag2}
(1 - \eta_j)  \mh( S_{inner, j} ) \leq \mathcal{H}^3 \Big\{ S_{inner,j} \cap
\mbox{exp} (1 \cdot \mathbf{D}_p (s) ) \Big\} s^{\epsilon}. 
\end{align}
In other words, most of the directions at 
$S_{inner, j}$, make it to $S_{inner, j+1}$ before 
hitting the cut locus. 

 Above we have shown that only 1 component
$\tilde{A}_j^{\circ}(3/2, s - 1/2 )$ has 
geodesics which make it to $S_{inner, j+1}$. 
Therefore all the radial geodesics in the other components 
must either hit the cut locus, or hit a different 
outer boundary component of $S_{outer, j}$. 
Let us call this set
\begin{align}
S_{bad, j} = 
\{ S_{inner, j} \setminus \mbox{exp} ( 1 \cdot \mathbf{D}_p (s) ) \}.
\end{align}
Therefore the other components are contained 
in the set 
\begin{align}
\tilde{A}_{bad, j} = \{ \gamma(t), 
1 \leq t \leq s : \gamma(t) \mbox{ is a radial 
geodesic with }  \gamma(1) \in S_{bad, j} \}.
\end{align}
 
By the coarea formula and area comparison, similar to (\ref{porl-4}) 
we have
\begin{align}
\label{gag1}
Vol (\tilde{A}_{bad, j}) \leq C \mathcal{H}^3
( S_{bad, j}) s^4. 
\end{align}
Next rewrite (\ref{gag2})
\begin{align*}
 \mh( S_{inner, j} ) & \leq \mathcal{H}^3 \Big\{ S_{inner,j} \cap
\mbox{exp} (s \cdot \mathbf{D}_p (s) ) \Big\} \frac{s^{\epsilon}}{(1-\eta_j)}\\
&= \big(  \mh (S_{inner,j}) - \mh (S_{bad,j}) \big) 
\frac{s^{\epsilon}}{(1 - \eta_j)},
\end{align*}
which implies
\begin{align*}
\mh (S_{bad,j}) \leq \mh( S_{inner,j}) \Big( 1 - \frac{1 - \eta_j}
{s^{\epsilon}} \Big). 
\end{align*}
Substituting into (\ref{gag1}), we have 
\begin{align*}
Vol (  (\tilde{A}_{bad, j}) \leq C \mh( S_{inner, j}).
\end{align*}
From the coarea formula, arguing as in (\ref{porl-3}), we have 
\begin{align*}
Vol ( \tilde{A}_j^{\circ}( 3/2, s-1/2) \geq
C \cdot \mathcal{H}^3 ( S(t_0) \cap \tilde{A}_j^{\circ}( 3/2, s-1/2)) (s-1), 
\end{align*}
where $t_0$ is the minimal area sphere (when rescaled to 
unit size) and $3/2 \leq t_0 \leq s-1/2$.  
Also, from comparison
\begin{align*}
 \mathcal{H}^3 ( S_{inner, j+1}) \leq \mathcal{H}^3 ( S(t_0) 
\cap \tilde{A}_j^{\circ}( 3/2, s-1/2)) ) s^{3+ \epsilon},
\end{align*}
since $S_{inner, j+1}$ must be covered by endpoints 
of radial geodesics going through 
$S(t_0) \cap \tilde{A}_j^{\circ}( 3/2, s-1/2)$ (geodesics hitting 
other components do not make it to $S_{inner, j+1}$). 

 Combining these inequalities, and using that 
our annulus satisfies $\mh (S_{inner,j+1}) \geq (1 - \eta_j)
\mh(S_{inner,j}) s^3$,
we have 
\begin{align*}
Vol ( \tilde{A}_j^{\circ}( 3/2, s-1/2) \geq 
&C \mathcal{H}^3 ( S(t_0) \cap \tilde{A}_j^{\circ}( 3/2, s-1/2)) (s-1)\\
& \geq C  \mathcal{H}^3 ( S_{inner, j+1})\\
& \geq C (1- \eta_j) \mathcal{H}^3 ( S_{inner, j}) s^3\\
& \geq C Vol (  (\tilde{A}_{bad, j}) \\
& \geq C Vol (  \tilde{A}_j( 3/2, s - 1/2 )
\setminus  \tilde{A}_j^{\circ}( 3/2, s - 1/2)). 
\end{align*}
\end{proof}

 We will now just work with the component $\tilde{A}_j^{\circ}( 3/2, s-1/2)$  
\begin{lemma}
\label{gongok}
We have
\begin{align}
Vol ( \tilde{A}_j(1,s) \setminus \tilde{A}_j^{\circ}( 3/2, s - 1/2 ))
\leq 3C \cdot Vol ( \tilde{A}_j^{\circ}( 3/2, s - 1/2 )). 
\end{align}
\end{lemma}
\begin{proof}
We have the inclusion
\begin{align*}
  \{ \tilde{A}_j(1,s)& \setminus \tilde{A}_j^{\circ}(3/2, s-1/2) \}\\
&\subset
 \{ \tilde{A}_j (3/2, s-1/2) \setminus \tilde{A}_j^{\circ}(3/2, s-1/2)\}
\cup
\{ \tilde{A}_j(1,s) \setminus \tilde{A}_j (3/2, s-1/2)\}.
\end{align*}
So from the previous 2 lemmas, we have 
\begin{align*}
 Vol (  \{ \tilde{A}_j(1,s)& \setminus \tilde{A}_j^{\circ}(3/2, s-1/2) \})\\
& \leq C  Vol ( \{  \tilde{A}_j^{\circ}(3/2, s-1/2)\})
+ C Vol(  \{ \tilde{A}_j (3/2, s-1/2)\} )\\
& \leq  3 C \cdot Vol ( \{  \tilde{A}_j^{\circ}(3/2, s-1/2)\}).
\end{align*}
\end{proof}

For $ \tilde{A}_j^{\circ}(3/2, s - 1/2)$, 
chose a maximal $1/4$-separated set, 
that is, choose points $p_{j,l} \in  
\tilde{A}_j^{\circ}(3/2, s - 1/2) , l = 1 \dots Q_j$ such that 
\begin{align*}
B(p_{j,l}, 1/8) \cap B(p_{j,l'}, 1/8) = \emptyset, 
\mbox{ for } l' \neq l,
\end{align*}
and such that 
\begin{align}
\label{gonginc}
\tilde{A}_j^{\circ} (3/2, s-1/2) \subset \cup_{l=1}^{Q_j}B(p_{j,l}, 1/4).
\end{align}
(For simplicity, we are assuming $s$ is very large, and 
we take $1/4$-separated sets. In general, we can take the 
separation to be some small multiple of $s$.)

We define the set $\tilde{K}_j =  \cup_{l=1}^{Q_j}B(p_{j,l}, 1/4)$.
From (\ref{gonginc}), and Lemma (\ref{gongok}), we have
\begin{align}
\label{gongoo}
Vol ( \tilde{A}_j(1,s) \setminus \tilde{K}_j )
\leq 3C \cdot Vol ( \tilde{K}_j ). 
\end{align}

For $i$ sufficiently large, from the curvature 
decay condition (\ref{decay12}), the curvature on 
$ \tilde{A}_j$ will be arbitrarily small. 
Using Bishop's volume comparison theorem, this means there 
exists a constant $C_2$ so that
\begin{align}
Vol( B(x,t)) \leq C_2 t^4,
\end{align}
for any $x \in \tilde{A}_j$, and $t < 1/2$.
Using the volume growth assumption (\ref{cond4}), 
we therefore have constants $C_1$ and $C_2$
so that (for $r$ sufficiently large)
\begin{align}
\label{nicevol}
C_1 t^4 \leq Vol( B(x,t)) \leq C_2 t^4.
\end{align}
This implies that $B(p_{j,l}, 1/4)$ can only hit a
uniformly bounded number of neighbors   $B(p_{j,l'}, 1/4)$,
since if 
\begin{align}
B(p_{j,l}, 1/4) \cap B(p_{j,l'}, 1/4) \neq \emptyset, 
\end{align}
then $B(p_{j,l'}, 1/4) \subset B(p_{j,l}, 1/2)$
and the latter must contain a bounded number 
of $1/4$-balls by (\ref{nicevol}).

 So we define graph $G_j$ with vertices 
$p_{j,l}$, $l = 1 \dots Q_j$, and $p_{j,l}$
is connected to $p_{j,l'}$ only if
$B(p_{j,l}, 1/4) \cap B(p_{j,l'}, 1/4) \neq \emptyset$.
By the above observation $G_j$ has a uniformly 
bounded number of edges at each vertex. 
Define a distance on the graph $d(p_l, p_{l'})
=$ minimal number of edges to traverse from $p_l$
to $p_{l'}$. 

 We begin to choose our sets, start with $p_1$, and 
choose the minimal integer $I$ such that 
\begin{align}
 I * Vol ( B(p_{1}, 1/4)) > 10 \pi^2 ( s^4 -1).
\end{align}
We then consider the union $ \tilde{K}_{j,1} = \cup_{l} B(p_{j,l}, 1/4)$ 
where the union is taken over all points $p_{j,l}$ such 
that $ d( p_{j,1}, p_{j,l}) \leq I$.
That is, we add all points $p_{j,l}$ at a bounded 
graph distance to $p_{j,1}$ until the union of the 
$1/4$-balls centered at these points has large enough 
volume. From the volume bounds (\ref{nicevol}),
$Vol ( \tilde{K}_{j,1})$ is uniformly bounded, and 
(iii) of the Proposition is satisfied.  

If there is no point $p_l$ such that $d(p_1, p_l) > 2I$
then we stop. 
All remaining points are within bounded 
graph distance, so the remaining portion 
has comparable volume, that is we must have 
\begin{align}
Vol ( \tilde{K}_j \setminus \tilde{K}_{j,1} ) \leq C Vol(\tilde{K}_{j,1}).
\end{align}

 Otherwise, we choose $p_2$ 
such that $d(p_1, p_2) > 2I$,
and repeat the process, our second set will 
be $ \tilde{K}_{j,2} = \cup_{l} B(p_{j,l}, 1/4)$ where 
the union is taken over all points $p_{j,l}$ such 
that $ d( p_{j,2}, p_{j,l}) \leq I$.

 Next, if there is no point  $p_{j,l}$ such that $d(p_{j,1}, p_{j,l}) > 2I$
and $d(p_{j,2}, p_{j,l}) > 2I$, then all the remaining 
balls are within bounded graph distance of 
$p_{j,1}$ and $p_{j,2}$. So what is remaining has 
small relative volume, and we must have the inequality 
\begin{align}
Vol ( \tilde{K}_j \setminus \{ \tilde{K}_{j,1}  \cup \tilde{K}_{j,2} \} ) 
\leq C Vol( \{\tilde{K}_{j,1}   \cup \tilde{K}_{j,2} \} ).
\end{align}

 We continue this procedure until it stops, 
and we are left 
with disjoint sets
\begin{align*}
\tilde{K}_{j,1}, \dots, \tilde{K}_{j,N_j} \subset 
\tilde{K}_j, 
\end{align*}
and such that we have 
\begin{align}
\label{gongop}
Vol ( \tilde{K}_j \setminus \{ \cup_{l=1}^{N_j} \tilde{K}_{j,l} \} ) 
\leq C Vol(  \{ \cup_{l=1}^{N_j} \tilde{K}_{j,l} \}  ).
\end{align}
We have the inclusion
\begin{align*}
  \{ \tilde{A}_j(1,s) \setminus  \{ \cup_{l=1}^{N_j} \tilde{K}_{j,l} \}
 \subset
 \{ \tilde{A}_j (1,s) \setminus \tilde{K}_j\}
\cup
\{ \tilde{K}_j \setminus   \{ \cup_{l=1}^{N_j} \tilde{K}_{j,l} \}  \}.
\end{align*}
From the inequalities (\ref{gongoo}) and (\ref{gongop}), this implies 
\begin{align}
\label{gongmain}
Vol ( \tilde{A}_j(1,s) \setminus \{ \cup_{l=1}^{N_j} \tilde{K}_{j,l} \} ) 
\leq 5C \cdot Vol(  \{ \cup_{l=1}^{N_j} \tilde{K}_{j,l} \} ),
\end{align}
which is (iv). 
\end{proof}

\begin{lemma}
\label{comparable}
Assume that
\begin{align}
\frac{ a_1 + a_2}{b_1 + b_2} < \delta,
\end{align}
with $a_1, a_2, b_1,$ and $b_2$ nonnegative. 
If $b_2 < C b_1$ then
\begin{align}
\frac{a_1}{b_1} < (1+C) \delta.
\end{align}
\end{lemma}
\begin{proof}
Since $a_2 > 0$, we have 
$\frac{ a_1}{b_1 + b_2} < \delta$, then 
$\frac{a_1}{ (1 + C) b_1} <  \delta$. 
\end{proof}
From Proposition \ref{l1con}, there exists a sequence 
$\delta_j \rightarrow 0$ as $j \rightarrow \infty$ such that
\begin{align*}
\frac{1}{ Vol( \tilde{A}_j)}
\Big(
\int_{\tilde{A}_j} \Big| \Delta (\tilde{\rho}_j)^2  - 8  \Big|
\tilde{dV_j} + 2 E_j \Big) < \delta_j.
\end{align*}
To simplify notation, let $ \bar{K}_j =  \cup_{l=1}^{N_j} \tilde{K}_{j,l}$.
We decompose $\tilde{A}_j(1,s)$ into disjoint sets 
\begin{align*}
\tilde{A}_j(1,s) =  \{ 
\tilde{A}_j(1,s) \setminus  \bar{K}_j \} \cup  \bar{K}_j.
\end{align*}
Using inequality (\ref{gongmain}), and 
applying Lemma \ref{comparable}, we have that 
\begin{align}
\label{comparable2}
\frac{1}{ Vol( \bar{K}_j ) }
\Big( 
\int_{ \bar{K}_j} \Big| \Delta (\tilde{\rho}_j)^2  - 8 \Big|
\tilde{dV_j} + 2 E_j \Big)< (1 + 5C) \delta_j.
\end{align}
\begin{proposition}
\label{l1con4}
For each $j$, there exists $l'$ (depending on $j$) so that on  
$K_{j,l'} \subset \tilde{A}_j$, we have  
\begin{align*}
\int_{\tilde{K}_{j,l'}} \Big| \Delta (\tilde{\rho}_j)^2 - 8 
\Big| \tilde{dV_j} + 2 E_{j,l'}  \rightarrow 0,
\end{align*}
as  $j \rightarrow \infty$, 
where $E_{j,l}$ denotes the error term, but only
integrated over $\mathcal{C} \cap \tilde{K}_{j,l'}$.
\end{proposition}
\begin{proof}
We write (\ref{comparable2}) as 
\begin{align*}
\frac{\int_{\bar{K}_j } \Big| \Delta (\tilde{\rho}_j)^2 - 8
\Big| \tilde{dV_j} + 2 E_j }{ Vol ( \bar{K}_j ) }
= \frac{ \sum_l \Big(
\int_{\tilde{K}_{j,l}} \Big| \Delta  (\tilde{\rho}_j)^2 - 8 \Big| \tilde{dV_j} 
+2 E_{j,l} \Big)}{ \sum_l Vol ( \tilde{K}_{j,l} ) } < \delta_j',
\end{align*}
since the error term is clearly additive under disjoint unions. 
This implies that 
\begin{align*} 
 \underset{l}{\mbox{min}} 
 \frac{  
\int_{\tilde{K}_{j,l}} \Big| (\Delta \tilde{\rho}_j)^2 - 8
\Big|\tilde{dV_j} + 2 E_{j,l}}
{ Vol ( \tilde{K}_{j,l} ) } < \delta_j'
\end{align*}
therefore for some $l'$, 
\begin{align*}
\int_{\tilde{K}_{j,l'}} \Big| \Delta (\tilde{\rho}_j)^2 - 8
\Big| \tilde{dV_j} + 2 E_{j,l'} < \delta_j'  Vol ( \tilde{K}_{j,l'} ).
\end{align*}
From (iii) of Proposition \ref{kchoose}, $ Vol ( \tilde{K}_{j,l'} )$ 
is uniformly bounded, and the proposition follows. 
\end{proof}
To simplify notation, we will now write $\tilde{K}_{j,l'}$ as just 
$\tilde{K}_j$, and $E_{j,l'}$ as $E_{j}$.
\begin{lemma}[Cheeger--Gromov compactness]
\label{Cheeger}
Let $\tilde{K}_j^{\circ}$ denote the interior of $\tilde{K}_j$. 
Then $(\tilde{K}_j^{\circ}, g/s^{2i})$ has a subsequence 
which converges uniformly in the $C^{1, \alpha}$ topology
on compact subsets to a flat Riemannian manifold $\tilde{K}_{\infty}^{\circ}$.
\end{lemma}
\begin{proof}
From the assumption on the curvature decay, on $\tilde{K}_j$ 
we have $|Rm| \leq \epsilon_1 \rightarrow 0$  
as $j \rightarrow \infty$.  
The assumption (\ref{cond4}) implies that there exists 
a constant $C_1 > 0$ so that 
\begin{align}
\label{cond3}
Vol(B(q,s)) \geq C_1 s^4,
\end{align}
for any $q \in \tilde{K}_j^{\circ}$, with $B(q,s) \subset 
\tilde{K}_j^{\circ}$. Together with the curvature bounds,
by \cite{CGT}, this implies that there exists a constant $C_2$ so 
that for any $x \in \tilde{K}_j^{\circ}$, 
\begin{align}
\label{injrad1}
\mbox{ inj}(x) \geq C_2 \rho(x), 
\end{align}
where $\mbox{ inj}(x)$ denotes the 
injectivity radius at $x$. Furthermore, the particular way 
we have chosen the sets (as union of balls), there is no 
collapsing. More precisely, for any $ \epsilon > 0$, consider 
\begin{align*}
\tilde{K}_j^{\epsilon}
= \{ x \in \tilde{K}_j : \tilde{dist}(x, \partial \tilde{K}_j) 
> \epsilon  \}. 
\end{align*}
From the definition of $\tilde{K}_j$ as unions of balls, for 
$\epsilon$ sufficiently small, $\tilde{K}_j^{\epsilon}$ 
is nonempty, and also connected.
We may then apply a suitable version of the Cheeger-Gromov 
convergence theorem, see \cite{Anderson}, \cite{Tian}.  
\end{proof}

It follows there exist $C^{1, \alpha}$ diffeomorphisms 
$\Phi_j : \tilde{K}_{\infty}^{\circ} \rightarrow \tilde{K}_j^{\circ}$ such 
that the metrics $(\Phi_j)^{*} g_j$
converge in the $C^{1, \alpha}$ topology on compact subsets 
of $\tilde{K}_{\infty}^{\circ}$. 
By passing to a subsequence, we may assume that the 
rescaled distance functions $\tilde{\rho}_j = (\rho/s^j) \circ \Phi_j$
converge to a function 
$\rho_{\infty}$ in the $C^{\alpha}$ topology on compact 
subsets of $\tilde{K}_{\infty}$,
since the distance function is Lipschitz with
Lipschitz constant $1$. Proposition \ref{l1con4} then implies  
\begin{align}
\label{l1con1again}
\int_{\tilde{K}_{\infty}} \Big| \Delta_j (\tilde{\rho}_j)^2 - 8 
\Big| \tilde{dV_j} + 2 E_j  \rightarrow 0,
\end{align}
as $j \rightarrow \infty$,
where $\tilde{\rho}_j = (\rho/s^j) \circ \Phi_j$ is the rescaled distance function 
on $\tilde{K}_j$, and $\Delta_j$ is the Laplacian 
with respect to $\tilde{g}_j = \Phi_j^{\ast}( g/s^{2j})$.
Similarly, we pull-back the error term $E_j$ to 
$\tilde{K}_{\infty}$ under the diffeomorphism $\Phi_j$. 
\begin{proposition}
\label{weaksol}
On $\tilde{K}_{\infty}^{\circ}$, 
$\rho_{\infty}$ is a weak solution of the equation 
\begin{align}
\label{theyare0}
\Delta_{\infty} (\rho_{\infty})^2 = 8,
\end{align}
\end{proposition}
\begin{proof}
For a test function $\phi \in C^{\infty}_c(\tilde{K}_{\infty}^{\circ})$ 
we have 
\begin{align*}
\int_{\tilde{K}_{\infty}} \rho_{\infty}^2 (\Delta_{\infty} \phi) \ 
dV_{\infty} 
& = \int_{\tilde{K}_{\infty}} (\Delta_{\infty} \phi) (\lim_{j \rightarrow \infty} 
\tilde{\rho}_{j}^2 \tilde{dV_j})\\
& =\lim_{j \rightarrow \infty} \int_{\tilde{K}_{\infty}} \tilde{\rho}_{j}^2 
(\Delta_{j} \phi) \ \tilde{dV_j},
\end{align*}
since $\tilde{g}_j \rightarrow g_{\infty}$ in $C^{1,\alpha}$
as $j \rightarrow \infty$. 
 Using the formula (\ref{intparts3}), we integrate by parts:
\begin{align*}
\int_{\tilde{K}_{\infty}} \rho_{\infty}^2 (\Delta_{\infty} \phi) \ 
dV_{\infty} 
& = \lim_{j \rightarrow \infty}\Big(  \int_{\tilde{K}_{\infty}} 
(\Delta_{j} \tilde{\rho}_{j}^2 ) \phi \ \tilde{dV_j} + 2 E_j \Big),\\
& =\lim_{j \rightarrow \infty} \Big( \int_{\tilde{K}_{\infty}} 
(\Delta_j \tilde{\rho}_{j}^2 - 8 ) 
\phi  \ \tilde{dV_j} + 2E_j \Big) + \int_{\tilde{K}_{\infty}} 8 \phi  \ 
dV_{\infty},
\end{align*}
where
\begin{align*}
E_j = \int_{\mathcal{C}_j \cap \tilde{K}_{\infty}} 
 \phi \Big(  \langle \tilde{\rho}_j (\nabla \tilde{\rho}_j)^{+}, 
\partial_{\mathcal{C}_j}^{+} \rangle 
+ \langle \tilde{\rho}_j (\nabla \tilde{\rho}_j)^{-}, \partial_{\mathcal{C}_j}^{-} 
\rangle \Big) d \mathcal{H}^3.
\end{align*}
From (\ref{l1con1again}), we have 
\begin{align*}
\lim_{j \rightarrow \infty} \Bigg| \int_{\tilde{K}_{\infty}} (\Delta_j
\tilde{\rho}_{j}^2 
- 8 ) 
\phi  \ \tilde{dV_j} \Bigg| 
\leq C \mbox{ sup}|\phi| \lim_{j \rightarrow 
\infty} \Vert \Delta_j \tilde{\rho}_{j}^2 
- 8 \Vert_{L^1(\tilde{K}_{\infty})}
= 0.
\end{align*}
Also, the error term above $E_j$ can be estimated in
the same manner by (\ref{errorterm}), so the error term 
also vanishes in the limit. Consequently, 
as $j \rightarrow \infty$, we have 
\begin{align*}
\int_{\tilde{K}_{\infty}} \rho_{\infty}^2 (\Delta_{\infty} \phi) \  
dV_{\infty}
=  \int_{\tilde{K}_{\infty}} 8 \phi  \ dV_{\infty}
\end{align*}
for any test function $\phi \in C^{\infty}_c(\tilde{K}_{\infty}^{\circ})$.
\end{proof}
In particular, from elliptic regularity, we 
conclude that $\rho_{\infty}$ is smooth in $\tilde{K}_{\infty}^{\circ}$, 
and the $\rho_{\infty}$ is a smooth solution of (\ref{theyare0}).

\begin{proposition}
\label{l1con3}
For the subsequence $\{ j \} \subset \{i\}$, 
as $j \rightarrow \infty$, 
\begin{align}
\int_{\tilde{K}_j} 
\Big| \Delta \tilde{\rho}_j  - \frac{3}{\tilde{\rho}_j} \Big| 
\tilde{dV_j} + E_j' \rightarrow 0,
\end{align}
where 
\begin{align*}
E_j' = \int_{\mathcal{C}_j \cap \tilde{K}_{j}} 
 \Big(  (\nabla \tilde{\rho}_j)^{+}, \partial_{\mathcal{C}_j}^{+} 
\rangle +  (\nabla \tilde{\rho}_j)^{-}, \partial_{\mathcal{C}_j}^{-} 
\rangle \Big) d \mathcal{H}^3.
\end{align*}
\end{proposition}
\begin{proof}
At any point where the distance function is smooth we have,  
\begin{align*}
\Delta (\rho^2) = 2 \rho \Delta \rho + 2 |\nabla \rho|^2, 
\end{align*}
Since $|\nabla \rho|^2 = 1$ almost everywhere, we have that
the equation
\begin{align*}
\Delta \rho^2 = 2 \rho \Delta \rho + 2, 
\end{align*}
holds almost everywhere. We then have 
\begin{align*}
\int_{\tilde{K}_j} 
 \Big| \Delta \tilde{\rho}_j  - \frac{3}{\tilde{\rho}_j} \Big| \tilde{dV_j}
&= \int_{\tilde{K}_j} 
\frac{1}{ 2 \tilde{\rho}_j} \Big| 2 \tilde{\rho}_j \Delta \tilde{\rho}_j - 6 \Big|
\tilde{dV_j}\\
& =  \int_{\tilde{K}_j} 
\frac{1}{ 2 \tilde{\rho}_j} \Big| \Delta (\tilde{\rho}_j)^2 - 8 \Big|
\tilde{dV_j}\\
& \leq \frac{1}{ 2}\int_{\tilde{K}_j} 
\Big| \Delta (\tilde{\rho}_j)^2 - 8 \Big| \tilde{dV_j},
\end{align*}
and the right hand side goes to zero as $j \rightarrow \infty$ 
by Proposition \ref{l1con4}. 
Clearly $E_j'$ is dominated by $E_j$ from 
Proposition \ref{l1con4}. 
\end{proof}
\begin{proposition}
On $\tilde{K}_{\infty}^{\circ}$, 
$\rho_{\infty}$ is a strong solution of the equation 
\begin{align}
\label{theyare1}
 \Delta_{\infty} \rho_{\infty} = \frac{3}{\rho_{\infty}}. 
\end{align}
\end{proposition}
\begin{proof}
Using Proposition \ref{l1con3}, 
the proof is identical to the proof of Proposition \ref{weaksol}.
\end{proof}

\begin{proposition}
On $\tilde{K}_{\infty}^{\circ}$, 
$\rho_{\infty}$ is a strong solution of the equation 
\begin{align}
\label{theyare2}
\nabla^2_{\infty} \rho_{\infty} 
= \frac{1}{\rho_{\infty}}( g^{\infty} - d \rho_{\infty} \otimes  d \rho_{\infty}). 
\end{align}
\end{proposition}
\begin{proof}
By (\ref{theyare0}) and (\ref{theyare1}), we have 
\begin{align*}
8 = \Delta (\rho_{\infty})^2 &= 2 ( \rho_{\infty} \Delta \rho_{\infty}
+ |\nabla \rho_{\infty}|^2)\\
& = 2 ( 3 +  |\nabla \rho_{\infty}|^2),
\end{align*}
therefore 
\begin{align}
\label{distfcn}
|\nabla \rho_{\infty}| \equiv 1. 
\end{align}
From the Bochner formula, since $\tilde{K}_{\infty}$ is 
flat, we have 
\begin{align*}
0 = \Delta | \nabla \rho_{\infty}|^2 &= ( \nabla \rho_{\infty} , 
\nabla \Delta \rho_{\infty})
+ | \nabla^2 \rho_{\infty}|^2\\
& =  \Big( \nabla \rho_{\infty} , \nabla \frac{3}{\rho_{\infty}} \Big) 
+ | \nabla^2 \rho_{\infty}|^2\\
& =  - \frac{3}{\rho_{\infty}^2} 
\Big( \nabla \rho_{\infty} , \nabla \rho_{\infty} \Big) 
+ | \nabla^2 \rho_{\infty}|^2\\
& =  - \frac{3}{\rho_{\infty}^2} + | \nabla^2 \rho_{\infty}|^2\\
& =  \Big| \nabla^2 \rho_{\infty} -
 \frac{1}{\rho_{\infty}}( g^{\infty} - d \rho_{\infty} \otimes  
d \rho_{\infty}) \Big|^2,
\end{align*}
and we are done. 

An alternative proof of this proposition is as follows:
from (\ref{distfcn}), it follows that $\rho_{\infty}$ is a 
distance function. Since $\tilde{K}_{\infty}$ is flat, 
the Hessian comparison theorem implies
\begin{align*}
\nabla^2_{\infty} \rho_{\infty} 
\leq \frac{1}{\rho_{\infty}}( g^{\infty} - d \rho_{\infty} \otimes  d \rho_{\infty}). 
\end{align*}
Together with (\ref{theyare0}), this implies the proposition. 
\end{proof}
Since $\tilde{K}_{\infty}$ is flat, equation (\ref{theyare2}) and 
the Gauss equations imply that level sets of $\rho_{\infty}$
have constant positive sectional curvature, and therefore are 
locally isometric to portions of Euclidean spheres.  
It follows that $\tilde{K}_{\infty}$ is contained 
in a Euclidean annulus (or quotient thereof),
$1 \leq r \leq s$, and therefore $ Vol( \tilde{K}_{\infty}) \leq 
\frac{\pi^2}{2}(s^4-1)$. But from (iv) of Lemma \ref{kchoose},
we chose $\tilde{K}_j$ so that 
$Vol( \tilde{K}_j) >  \frac{\pi^2}{2}(s^4-1)+1$.
Since the convergence is smooth, we have 
$Vol( \tilde{K}_{\infty}) \geq \frac{\pi^2}{2}(s^4-1)  +1$, 
a contradiction. Thus there exists a 
constant $C$ so that $Vol( \tilde{A}_j(1,s)) \leq C$,
for each annulus in the subsequence chosen 
in Proposition \ref{annchoose}.
Moreover, we have proved that the distance function
converges to the Euclidean distance function 
on this subsequence, so not only is the 
volume bounded, the area of level sets is 
bounded. 
Therefore, for our original subsequence $s^j$, there exists a constant $C$ 
so that for all $j$ we have  
$ \mathcal{H}_{A_j}(s^j) \leq C s^{3j}$. We conclude that for all $i$, 
\begin{align}
\label{goodest}
\mathcal{H}_{A_i}(s^i) \leq C s^{3i},
\end{align} 
(since for $i$ not in
our subsequence, we will have $\mathcal{H}_{A_i}(s^{i+1}) \leq 
\mathcal{H}_{A_i}(s^i)s^3$).
Finally, given any $r$ sufficiently large, 
choose $i$ so that $$s^i \leq r \leq s^{i+1}.$$
From Lemma \ref{radialcomp} and (\ref{goodest})
and Lemma \ref{components}, we have 
\begin{align}
\mathcal{H}_{A_i}(r) \leq \mathcal{H}_{A_i}(s^i) 
s^{3+ \epsilon} \leq \mathcal{H}_{A_i}(s^i) s^3 (1 + \epsilon)
\leq C s^{3i} \leq Cr^3,
\end{align}
therefore by the coarea formula there exists a constant 
$C$ so that for $s^i \leq r < s^{i+1}$, 
\begin{align}
\label{vg2}
Vol(B(p, r) \cap A_i(s^i, s^{i+1}) ) \leq C r^4.
\end{align}
An examination of the proof shows that 
for $r$ large, we may take the constant 
$C$ in (\ref{vg2}) to be close to the 
corresponding Euclidean constant, that is 
\begin{align}
\label{goaa} 
Vol(A_i(s^i, s^{i+1})) \leq (\omega_4 + \alpha_i) ( s^{4(i+1)} - s^{4i}), 
\end{align}
where $\alpha_i \rightarrow 0$ as $i \rightarrow \infty$
(by Euclidean constant we mean $Vol(B(0,r)) = \omega_4r^4 = (\pi^2/2)r^4$ in 
$\mathbf{R}^4$). This is because for $j$ in the subsequence, our proof 
above shows that $\Delta \tilde{\rho}_j \rightarrow (3 / \tilde{\rho}_j)$ as
$j \rightarrow \infty$, so for any subsequence
$\{j'\} \subset \{ j \}$, there will be a 
subsequence $ \{k \} \subset \{j'\}$ so that 
$\tilde{A}_{k}( 1,s)$ will converge to a
Euclidean annulus or quotient thereof (note that when repeating 
our above $L^1$-convergence argument, we do not need to cut as in 
Proposition \ref{kchoose}, since we have already proved the estimate 
(\ref{vg2})). We conclude that for $j$ sufficiently large, 
\begin{align}
\label{ohgo}
\mathcal{H}_{A_j}(s^j) \leq (4 \omega_4 + \alpha_j') s^{3j}, 
\end{align}
where $\alpha_j' \rightarrow 0$ as $j \rightarrow \infty$.
Recall again that for $i$ not in the subsequence 
$\{j\}$, we have $\mathcal{H}_{A_i}(s^{i+1}) \leq 
\mathcal{H}_{A_i}(s^i)s^3$. This fact, inequality (\ref{ohgo}), 
area comparison and the coarea formula then imply (\ref{goaa}). 

We have derived a bound for the volume of the original sequence 
of annular components $A_j(s^j,s^{j+1})$, but we are free
to change the starting radius of the dyadic sequence. The above 
argument therefore shows that we have volume bounds on any annular components 
of fixed ratio of inner and outer radius, that is for any $1 \leq r_0 < s$, 
\begin{align}
\label{goab} 
Vol(A_i(r_0s^i, r_0 s^{i+1})) \leq (\omega_4 + \alpha_i) 
\Big( (r_0 s)^{4(i+1)} - (r_0 s)^{4i} \Big), 
\end{align}
where $\alpha_i \rightarrow 0$ as $i \rightarrow \infty$.
Since $s$ was arbitrary, the above will hold for all 
$s$, although the sequence $\alpha_i$ might depend on $s$.

We next consider the cone at infinity in the 
direction of our chosen sequence of annuli. 
For any $n > 0$, and $t >0$, 
let $(\tilde{A}^t(1/n, n), \tilde{g}) = 
(A(t/n, nt), g/r^2)$ denote the component 
of the annulus $\{x : t/n < \rho(x) < nt\}$,
which touches some annular component in our 
chosen sequence, with the rescaled metric $g/t^2$. 
From (\ref{decay1}) we have 
\begin{align}
\label{decay111}
\underset{A^t}{\mbox{sup}}\tilde{|Rm|} \leq \epsilon_3(t),
\end{align} 
where $\epsilon_3(t) \rightarrow 0$ as $t \rightarrow \infty$.
As in Lemma \ref{Cheeger},  it follows that some subsequence 
$(\tilde{A}^{t_i}, \tilde{g}^{t_i})$ converges 
in the $C^{1, \alpha}$ topology on 
compact subsets to a flat limit space $(A^{\infty}(1/n,n), g_{\infty})$.
By taking subsequences, we construct limit spaces with natural inclusions 
$ A^{\infty}(1/m,m) \subset A^{\infty}(1/n,n)$
for $m < n$. We may then take a subsequence 
$n_i$ and obtain in the limit a flat 
connected Riemannian manifold $A^{\infty}(0, \infty)$ 
with $A^{\infty}(1/n,n) \subset A^{\infty}(0, \infty)$
for any $n$. Note there is a distance function 
$\rho_{\infty}$ obtained as a limit of the (rescaled) original 
distance function $\rho(p, \cdot)$. From (\ref{goab}), follows the estimate
\begin{align}Vol( A_{\infty}(r_1,r_2))
\label{vg3}
\leq \omega_4 (r_2^4 - r_1^4),
\end{align}
for any $0 < r_1 < r_2.$

We claim that level sets $S_{\infty}(s) \equiv \{ \rho_{\infty} = s \} \subset 
A^{\infty}(0, \infty)$ are connected. 
First, $ A^{\infty}(0, \infty)$ must be connected, since we 
have constructed this tangent cone at infinity by taking
sequences of connected components of annuli which touch, 
that is, these sequences correspond to an end of the manifold.   
If $S_{\infty}(s)$ was disconnected for $s$ small, then we would find 
a connected annulus in  $ A^{\infty}(0, \infty)$
with more than one component in the beginning portion, or 
a ``bad'' annulus in the above terminology. 
This would correspond to an infinite sequence of 
bad annuli in $X$, which contradicts the 
Betti number estimate $b_1(X) < \infty$.

Since $S_{\infty}(s)$ is connected for $s$ small,
(\ref{vg3}) and (\ref{cond4}) then imply an estimate 
\begin{align}
\label{dmest}
\mbox{diam} \{ S_{\infty}(s) \} \leq C s.
\end{align}
A standard analysis of the developing map for flat structures
(see \cite{EschenburgSchroeder}, \cite{Tian}, \cite{Anderson}), 
implies that we may add a point to $A^{\infty}(0,\infty)$ to obtain a complete 
length space $\overline{A^{\infty}}$, 
and we conclude that $\overline{A^{\infty}}$ is isometric to a cone 
on $S^{3}/ \Gamma$ where $\Gamma \subset SO(4)$ is a finite 
subgroup of isometries of $S^3$. 

 An alternative proof of this step, not relying on the Betti
number estimate, is as follows. The estimate (\ref{dmest})
will hold for each component of $S_{\infty}(s)$.
Using an analysis of the developing map, and the estimate (\ref{vg3}), the
universal cover of $A_{\infty}(0,\infty)$ must be 
$ \mathbf{R}^n \setminus \{p_0, \dots, p_k\}$,
with distance function being the distance to this 
finite set, that is, $\rho_{\infty}(\cdot) = \min_i d(p_i,\cdot)$
(see \cite{EschenburgSchroeder}).
This implies that $\overline{A^{\infty}}$ is   
a flat cone $\mathcal{C}(S^3/\Gamma)$ with finitely 
many points $\{p_1, \cdots, p_k\}$ identified to the vertex $p_0$, 
and distance function  $\rho_{\infty}(\cdot) = \min_i d(p_i,\cdot)$. 
The volume growth estimate (\ref{vg3}) then implies that, 
$\rho_{\infty}(\cdot) = d(p_0,\cdot) $, therefore 
$\overline{A^{\infty}}$ is isometric to the flat cone $C(S^3 / \Gamma)$. 

We have shown that every tangent cone, in the 
direction of our initial sequence of annuli, 
is isometric $C(S^3 / \Gamma)$. It is standard to
conclude that $ \Gamma$ is unique and to construct
an ALE coordinate system of order zero,
we refer the reader to \cite{Tian}, \cite{Anderson}
for details.

 This implies that our sequence of annuli must eventually 
meet a Euclidean end of the 
manifold. Since that the initial sequence of annuli can 
be chosen arbitrarily, there must be finitely many ends. 
That is, assuming there are infinitely many ends, then we can 
choose a sequence of components of annuli which touch, 
and for which there are always infinitely many ends past any given annular 
component in the sequence.
This is a contradiction, since such a sequence 
could not lead to a Euclidean end. 
This finishes the proof of Theorem \ref{bigthm}. 
\section{Kato inequalities and Ricci decay}
\label{proof2section}
In this section, we discuss some improved Kato inequalities, 
which we use to improve Ricci curvature decay.
\begin{lemma}
If $(X,g)$ satisfies  $\delta W^+ = 0$ and has 
zero scalar curvature, then 
\begin{align}
\label{rickato}
\big| \nabla |Ric| \big|^2 \leq \frac{2}{3} | \nabla Ric | ^2,
\end{align}
at any point where $|Ric| \neq 0$. 
\end{lemma}
\noindent
{\bf{Remark.}}
This is rather surprising, because if one assumes moreover 
that $\delta W = 0$, the best constant 
is still $2/3$. In an earlier version of this paper, 
we had an improved Kato inequality in this case, 
but our constant was not optimal. 
The authors are indebted to Tom Branson for 
showing us how the best constant in this case
follows easily from his general theory 
of Kato constants in \cite{Branson}. 

\begin{proof}
\def\tfs{{\rm TFS}}
\let\a=\alpha
\let\b=\beta
\let\f=\varphi
\let\F=\Phi
\let\l=\lambda
\let\L=\Lambda
\let\m=\mu
\let\N=\nabla
\let\s=\sigma
\let\y=\psi
\let\Y=\Psi
\def\tc{\tilde c}
\def\tb{\tilde\beta}
\def\dfrac#1#2{\displaystyle{\frac{#1}{#2}}}

We assume familiarity with the paper of Branson \cite{Branson}.
We are in dimension 4, so weight labels have only 2 entries. 
Since $Ric$ is traceless, $Ric$ is a section of 
the $(2,0)$ bundle, which we denote by $TFS^2$. 
The covariant derivative $\nabla Ric$ is then a section of 
$T^*X \otimes TFS^2$. This decomposes into irreducible pieces 
under the action of the orthogonal group, and
the selection rule targets are $(3,0),\;(2,1),\;(2,-1),\;(1,0)$.  
These correspond to trace-free symmetric $3$-tensors, 
the self-dual and anti-self-dual parts of
Codazzi tensors, and $\L^1$.  The assumption is that $Ric$ is 
in the kernel of the gradients targeted at $(1,0)$ and $(2,1)$. 

A calculation shows that
$$
\begin{array}{l}
s_{(2,1)}=s_{(2,-1)}=-\frac12\,, \\
s_{(3,0)}=-\frac72\,,\ \ s_{(1,0)}=\frac52\,.
\end{array}
$$
We next compute some of the $\tilde c$ numbers.  Note that $t(\l)=2$ for 
$\l=(2,0)$; substituting and using the $s$ numbers above,
$$
\begin{array}{l}
\tc_{(2,1)}=\tc_{(2,-1)}=\frac1{18}\,, \\
\tc_{(1,0)}=-\frac1{18}\,.
\end{array}
$$

The Kato constant is expressed as a minimum over the branching products of
$(2,0)$.  These branching products are the SO$(3)$-modules $(2)$,
$(1)$, and $(0)$, which have dimensions $5$, $3$, and $1$, respectively.
The minimum works out as
$$
\min_{\b=0,1,2}\left(
\tc_{(2,1)}
\left(\tb^2-\left(\frac12\right)^2\right)
+\tc_{(1,0)}\left(\tb^2-\left(\frac52\right)^2\right)
\right)
$$
where $\tb=\b+\frac12$. 
This is 
$$
\min\left(\frac13,\frac13,\frac13\right)=\frac13.
$$
Therefore the Kato constant is $1/3$, which is 
\begin{align*}
\big| \nabla |Ric| \big|^2 \leq (1 - \frac{1}{3}) | \nabla Ric | ^2
= \frac{2}{3} | \nabla Ric | ^2.
\end{align*}

We remark that if one takes the stronger assumption that $\delta W=0$, 
the one must add in $(2,-1)$ terms above, equal to the $(2,1)$ terms.
The result is
$$
\min\left(\frac13,\frac49,\frac23\right)=\frac13,
$$
which is still the same Kato constant.

\end{proof} 
\begin{remark} We would like to thank the referees for pointing out that 
the same constant follows from the methods in \cite{CGH}.
The case considered in Lemma \ref{rickato} is
exactly the case $r=s=2$ in the last line of the table  
on \cite[page 253]{CGH}, giving immediately the best constant $2/3$.
\end{remark}
\begin{proposition} 
Under the assumptions of Theorem \ref{decayale},
in case (a) we have
\begin{align}
\underset{S(r)}{sup } \ |Ric| = O(r^{- \alpha}), \mbox{ for any } \alpha < 4.
\end{align}
\end{proposition}
\begin{proof} 
From (\ref{generaleqn}), there exists a constant $C$ so that 
\begin{align*}
|\Delta Ric | \leq C |Rm||Ric|.
\end{align*}
Noting that 
\begin{align*} 
\frac{1}{2} \Delta | Ric|^2 = \langle Ric, \Delta Ric \rangle
+ |DRic|^2 = |Ric| \Delta |Ric| + |d|Ric||^2,
\end{align*}
and using the Kato inequality (\ref{rickato}), we obtain
\begin{align*}
\langle Ric, \Delta Ric \rangle \leq  |Ric| \Delta |Ric|
- \frac{1}{3}|D Ric|^2
\end{align*}
Next we have 
\begin{align*}
-|Ric| | \Delta Ric| \leq - | \langle Ric, \Delta Ric \rangle |
\leq  \langle Ric, \Delta Ric \rangle 
\leq  |Ric| \Delta |Ric| - \frac{1}{3} |D Ric|^2,
\end{align*}
so that when $|Ric| > 0$, 
\begin{align}
\label{Riceqn2}
\Delta | Ric| \geq - C |Rm| |Ric| + \frac{1}{3} \frac{|D
  Ric|^2}{|Ric|}. 
\end{align}
We obtain
\begin{align*} 
\Delta |Ric|^{1/2} &= - \frac{1}{4}
|Ric|^{-3/2} | d |Ric||^2 
+ \frac{1}{2} |Ric|^{-1/2} \Delta | Ric|\\
&\geq \ - \frac{1}{6} 
|Ric|^{ -3/2} | D Ric|^2 + \frac{1}{2} |Ric|^{- 1/2} 
 \left( - C |Rm| |Ric| + \frac{1}{3} \frac{|D
  Ric|^2}{|Ric|} \right)\\ 
& \ = - \frac{1}{2} C |Ric|^{1/2} |Rm|. 
\end{align*}
Therefore $|Ric|$ satisfies weakly the differential inequality
\begin{align} 
\label{impriceqn}
\Delta |Ric|^{1/2} \geq - C |Ric|^{1/2} |Rm|.
\end{align}
Recall that $D(R) = X \setminus B(p,R)$.
\begin{lemma}
\label{goth}
There exist constants $C, R$, depending upon $C_S$ so that 
for $r > R$ we have 
\begin{align}
\label{mos2ann}
\left\{ \int_{D(2r)} |Rm|^4 dV_g \right\}^{1/2}
\leq
\frac{C}{r^2} \int_{D(r)} |Rm|^2 dV_g.
\end{align}
\end{lemma}
\begin{proof}
This is proved as in (\ref{mos2}) above, except now we choose
the cutoff function $\phi$ such that 
$ \phi \equiv 1$ on $D(2r) \setminus D(r')$, $\phi = 0$ on $B(r) \cup D(2r')$, 
and $| \nabla \phi| < C( r^{-1} + r'^{-1})$, and then let 
$r' \rightarrow \infty$. 
\end{proof}

From Lemma \ref{goth}, and equation (\ref{impriceqn}), 
we may apply \cite[Proposition 4.8 (1)]{BKN} to conclude $|Ric| = O(r^{- \alpha})$
for any $\alpha < 4$ in the self-dual or 
anti-self-dual case.
\end{proof}
  Next, we consider the case of harmonic curvature
in dimension 4. The following lemma is standard, 
since harmonic curvature is the Riemannian 
analogue of Yang-Mills (see \cite{Rade}). 
\begin{lemma}
If $(X,g)$ satisfies  $\delta Rm = 0$ and 
has zero scalar curvature, then 
\begin{align}
\label{riemkato}
\big| \nabla |Rm| \big|^2 \leq \frac{2}{3} | \nabla Rm | ^2,
\end{align}
at any point where $|Rm| \neq 0$. 
\end{lemma}

\begin{proposition} 
Under the assumptions of Theorem \ref{decayale},
in case (b) we have
\begin{align}
\underset{S(r)}{sup } \ |Rm| = O(r^{- \alpha}), \mbox{ for any } \alpha < 4.
\end{align}
\end{proposition}
\begin{proof} 

A computation similar to above shows we have the improved equation
\begin{align} 
\Delta |Rm|^{1/2} \geq - C |Rm|^{1/2} |Rm|. 
\end{align}
Applying the method in Section $4$ of 
\cite{BKN} we find that $|Rm| = O(r^{- \alpha})$
for any $\alpha < 4$.
\end{proof}
 
 This completes the proof of Theorem \ref{decayale} for the case of
harmonic curvature, since by \cite[Theorem 1.1]{BKN}, 
$(X,g)$ is ALE of order $\tau$ for any $\tau < 2$. 
\begin{remark}
We remark that a 
version of Theorem \ref{decayale} in the locally 
conformally flat case was considered in \cite{CarronHerzlich}
(however, in that work the volume growth estimate 
(\ref{vga_i}) was insufficiently justified).
\end{remark}
\section{Asymptotically locally Euclidean}
\label{removal}
In this section, we complete the proof of Theorems
\ref{decay} and \ref{decayale}.
We first have a lemma on lower volume growth
\begin{lemma}
\label{volbelow}
Let $(X,g)$ be complete, noncompact, and 
assume that the Sobolev inequality (\ref{siq}) is satisfied
for all $f \in C^{0,1}_c(X)$.  
Then there exists a constant $C_1 > 0$ 
depending only upon $C_S$ so that
\begin{align}
\label{precond4}
Vol(B(q,s)) \geq C_1 s^4,
\end{align}
for any $q \in X$, and all $s \geq 0$. 
\end{lemma}
\begin{proof}
Fix $r > 0$, and consider the function 
\begin{align*}
f(x) = \left\{ \begin{array}{ll}
d(q,x) - r &  d(q,x) \leq r\\
0 &  \mbox{otherwise}.
\end{array} \right.
\end{align*} 
Note $f$ is Lipschitz, with $|\nabla f| = 1$ almost 
everywhere in $B(q,r)$, so in particular $f \in C^{0,1}_c(X)$. 
From the Sobolev inequality (\ref{siq}), we obtain
\begin{align*}
\frac{r}{2} Vol(B(q,r/2))^{1/4} \leq Vol(B(q,r))^{1/2}.
\end{align*}
The lemma follows by iterating this inequality (see \cite[Lemma 3.2]{Hebey}). 
\end{proof}
\begin{remark} The connection between Sobolev constant and lower 
volume growth has appeared many times in the literature,
see for example \cite[Proposition 2.1]{Akutagawa},
\cite{Carron2}. 
\end{remark}
\begin{proposition} 
\label{ale0thm}
Let $(X,g)$ satisfy the assumptions 
of Theorem \ref{decay}. Then $(X,g)$ has 
finitely many ends, and each end of $(X,g)$ is ALE of 
order zero. 
\end{proposition}
\begin{proof}
The volume growth condition (\ref{cond4}) 
follows from Lemma \ref{volbelow}.
The curvature decay condition (\ref{decay1}) follows from 
Theorem \ref{higherlocalregthm2}. The proposition then follows 
from Theorem \ref{bigthm}. 
\end{proof}
Proposition \ref{ale0thm} completes the proof
of Theorem \ref{decay}.
For Theorem \ref{decayale}, we assume furthermore that 
$(X,g)$ is of type (a), (b), or (c). 
We begin by listing some of the properties of $(X,g)$ that we have proved:
$(X,g)$ is a complete noncompact 4-manifold with 
finitely many ends with base point $p$, satisfying
\begin{align}
\label{1}
&\mbox{Each end of $X$ is ALE of order zero,} \\
\label{2}
&\Vert Rm \Vert_{L^2(X)} < \infty, \\
\label{3}
&\underset{S(r)}{\sup} |Ric| = O(r^{-(2 +\alpha)}), \mbox{ for any } \alpha < 2,\\
\label{4}
&\underset{D(2r)}{sup}|Rm| \leq
\frac{C}{r^2} \left\{ \int_{D(r)} |Rm|^2 dV_g \right\}^{1/2},\\
\label{5}
&\int_{D(r)} | \nabla Ric|^2 dV_g 
\leq \frac{C}{r^{2 + 2 \alpha}}.
\end{align}

Note that (\ref{3}) was proved in Section \ref{proof2section}, 
(\ref{4}) is the 
$\epsilon$-regularity Theorem \ref{higherlocalregthm2}.
Note that (\ref{5}) is proved as in (\ref{moslater}), but 
now taking a cutoff function supported in the complement 
of a ball (see the proof of Lemma \ref{goth}). 

  We will show that conditions (\ref{1})-(\ref{5}) imply that 
$|Rm| = O(r^{-(2+\alpha)})$. 
For simplicity, assume $(X,g)$ has only one end  
(in general do the same argument for each end). 
The proof is very similar to the 
argument given in Section 4 of \cite{Tian}, and that argument was inspired 
by the proof of Uhlenbeck for Yang-Mills connections in 
\cite{Uhlenbeck2}. 
Let $\tilde{A}$ be the Levi-Civita 
connection form of the metric $g$, where the covariant 
derivative is given by $\tilde{D} = d + \tilde{A}$. 
Note that $\tilde{A} \in C^{1, \alpha}( X, so(4) \times R^4)$.
Since $g$ is ALE of order zero, we transform from Euclidean 
to spherical coordinates, and any one-form $A = (A_r, A_{\psi})$
splits into radial and spherical parts. 
The following lemma, which is proved in 
\cite[Lemma 4.1]{Tian}), shows we can choose a good gauge, and 
is essentially an application of the Implicit Function Theorem.
\begin{lemma}
\label{Tian4.1}
Let $r$ be sufficiently large. Then there is a gauge 
transformation $u$ in $C^{\infty}(A(r,2r), so(4))$ 
with the property that if 
$D = e^{-u} \cdot \tilde{D} \cdot e^{u} = d+ A$, then 
\begin{align}
D^{\ast}A &= 0 \mbox{ on } \overline{ A(r,2r)}\\
d^{\ast}_{\psi} A_{\psi} &= 0 \mbox{ on } \partial A(r,2r)\\
\int_{A(r,2r)} & A( \nabla_F r(x)) dV_g = 0,
\end{align}
where $d^{*}, d_{\psi}^{*}$ are the adjoint operators 
of the exterior differentials on $A(r,2r)$, $\partial A(r,2r)$
with respect to $g$, respectively, $\nabla_F$ denotes the 
standard gradient, and $dV_g$ is the volume form with respect to $g$. 
Moreover, we have 
\begin{align}
\underset{A(r,2r)}{sup} (\Vert A \Vert_g(x)) \leq \frac{\epsilon_2(r)}{r},
\end{align}
where $\epsilon_2(r)$ is a decreasing function of $r$ with 
$\lim_{r \rightarrow 0} \epsilon_2(r) = 0$. 
\end{lemma}
The next lemma shows how the connection form decays in 
the Hodge gauge, and the proof is exactly as in \cite[Lemma 4.2]{Tian}.
\begin{lemma}
\label{Tian4.2}
Let $A$ be the connection form given in Lemma \ref{Tian4.1}. 
For $r$ sufficiently large, we have
\begin{align}
\underset{A(r,2r)}{sup} \Vert A \Vert_g (x) 
&\leq C r \underset{A(r,2r)}{sup} \Vert R_{A} \Vert_g (x) \\
\int_{A(r,2r)}  \Vert A \Vert_g^2 (x)dV_g 
&\leq C r^2 \int_{A(r,2r)}  \Vert R_{A} \Vert_g^2 (x)dV_g. 
\end{align}
\end{lemma}
The following lemma gives an improvement on the 
decay of the full curvature tensor, given the 
Ricci decay. This was proved in \cite[Lemma 4.3]{Tian} 
for the Einstein case, but with some extra work 
we prove here that the only improved Ricci 
decay is necessary. 
\begin{lemma}
\label{bali}
There exists $\beta > 0$ such that for $r$ sufficiently 
large, we have 
\begin{align}
\label{pdc}
\underset{D(2r)}{sup} \Vert Rm_g \Vert_g
\leq \frac{C}{r^{2 + \beta}}.
\end{align}
\end{lemma}
\begin{proof}
Choose  $r_0$ large, and let $r_i = 2 r_{i-1}$.
Let $A_i$ be the connection on $A(r_{i-1}, r_{i})$ 
from Lemma \ref{Tian4.1}. Then 
\begin{align*}
d_{\psi}^* A_{i \psi}|_{ \partial A(r_{i-1}, r_{i})} =0, \\
d_{\psi}^* A_{ (i-1) \psi}|_{\partial A(r_{i-2}, r_{i-1})} = 0, 
\end{align*}
so the restrictions $A_{i \psi}$ and $A_{(i-1) \psi}$ 
differ by a constant gauge on $\partial B(r_{i-1})$,
and we may therefore assume that 
\begin{align*}
A_{i \psi}|_{ \partial B( r_{i-1})}= 
A_{ (i-1) \psi}|_{\partial B(r_{i-1})}.
\end{align*}
Letting $\Omega_i = A(r_{i-1}, r_i)$, we compute
\begin{align*}
\int_{\Omega_i} \Vert R_{A_i} \Vert^2_{g} dV_g
&= \int_{\Omega_i}  \langle d A_i + A_i \wedge A_i, R_{A_i}
\rangle_g dV_g
= \int_{\Omega_i} \langle D_i A_i - [ A_i, A_i], R_{A_i} \rangle_g
dV_g\\
& = - \int_{\Omega_i} \langle [ A_i, A_i], R_{A_i} \rangle_g
dV_g - \int_{\Omega_i} \langle  A_i, D_i^* R_{A_i} \rangle_g
dV_g\\
& \ \ \ \  - \int_{S_i} \langle A_{i \psi}, (R_{A_i})_{r \psi} \rangle_g
d \sigma_g + \int_{S_{i-1}}  \langle A_{i \psi}, (R_{A_i})_{r \psi} \rangle_g
d \sigma_g, 
\end{align*}
where $D_i = d + A_i$, $S_i = \partial B(r_i)$, and $d \sigma$ is the 
induced area form. 

Next we sum over $i$, 
and since $(R_{A_i})_{r \psi}= (R_{A_{i+1}})_{r \psi}$ on $S_i$, 
we obtain
\begin{align*}
\int_{D(r_0)} \Vert Rm_g \Vert^2_g dV_g 
&= \sum_{i=1}^{\infty} \int_{\Omega_i} \Vert Rm_g \Vert^2_g dV_g\\ 
&= \sum_{i=1}^{\infty} \int_{\Omega_i} \Vert R_{A_i} \Vert^2_g dV_g\\
& = - \sum_{i=1}^{\infty} \int_{\Omega_i} \langle 
[ A_i, A_i], R_{A_i} \rangle_g dV_g  
+  \int_{\partial B(r_0)} \langle A_{1 \psi}, (R_{A_1})_{r \psi} 
\rangle_g d \sigma_g \\
& \ \ \ \  - \sum_{i=1}^{\infty} \int_{\Omega_i} \langle  A_i, D_i^* R_{A_i} 
\rangle_g dV_g.
\end{align*}

Since $A_i$ is equivalent to $\tilde{A}$ by a gauge transformation, 
we have $D^* R_{A_i} = D^* R_{\tilde{A}}$. 
Using the inequality $ab \leq \delta r^{-2}a + \delta^{-1} r^{2}b$, 
assumption (\ref{5}),
Lemma \ref{Tian4.2}, and the Bianchi identity 
$D^* Rm = d^{\nabla} Ric$, we estimate the last term as 
\begin{align*}
\begin{split}
 \int_{\Omega_i} \langle  A_i, D_i^* R_{\tilde{A}} 
\rangle_g dV_g
& \leq   {\delta}{r_i^{-2}}\int_{\Omega_i} \Vert A_i 
\Vert_g^2 dV_g 
+  C  r_i^2 \int_{\Omega_i} \Vert \nabla Ric 
\Vert_g^2 dV_g \\
&  \leq   \delta \int_{\Omega_i} \Vert R_{A_i} 
\Vert_g^2 dV_g 
+  C  r_i^{-2 \alpha}.
\end{split}
\end{align*}
Summing this over $i$, we obtain
\begin{align*}
\sum_{i=1}^{\infty}  \int_{\Omega_i} \langle  A_i, D_i^* R_{\tilde{A}} 
\rangle_g dV_g
& \leq \delta \int_{D(r_0)} \Vert Rm_{g} \Vert_g^2 dV_g 
+ C \sum_{i=1}^{\infty}   r_i^{- 2 \alpha}\\
& \leq \delta \int_{D(r_0)} \Vert Rm_{g} \Vert_g^2 dV_g 
+ C r_0^{-2 \alpha}.
\end{align*}
Also we have 
\begin{align*}
\Big|  \sum_{i=1}^{\infty} \int_{\Omega_i} \langle 
[ A_i, A_i], R_{A_i} \rangle_g dV_g \Big| 
&\leq \sum_{i=1}^{\infty} C \underset{x \in \Omega_i}{sup} ( \Vert
R_{A_i} \Vert_g (x) )\int_{\Omega_i} \Vert A_i \Vert^2_g dV_g\\
& \leq  \sum_{i=1}^{\infty} C \epsilon(r_{i-1}) \int_{\Omega_i} 
\Vert R_{A_i} \Vert^2_g dV_g\\
& \leq C \epsilon(r_0) \int_{D(r_0)} \Vert Rm_{g} \Vert^2_g dV_g.
\end{align*}
Combining the above we obtain
\begin{align*}
(1 -  C \epsilon(r_0) - \delta)  \int_{D(r_0)} \Vert Rm_{g} \Vert^2_g dV_g
\leq  \int_{\partial B(r_0)} \langle A_{1 \psi}, (R_{A_1})_{r \psi} 
\rangle_g d \sigma_g +  C r_0^{- 2 \alpha}.
\end{align*}
We have shown in Section \ref{euclideanvolumegrowth} 
that $(\partial B(r), \frac{1}{r^2} g$) 
converges uniformly to $( S^3 / \Gamma, g_{0})$, a quotient of 
the unit sphere with the standard metric, in the 
Gromov-Hausdorff topology.
Therefore, by the proof of Corollary 2.6 
in \cite{Uhlenbeck2}, we may find a decreasing 
function $\epsilon'(r)$ of $r$ with 
$\lim_{r \rightarrow 0} \epsilon'(r) = 0$ such that 
\begin{align}
\int_{\partial B(r)} \Vert A_{1\psi} \Vert^2_g dV_g
\leq (2 - \epsilon'(r))^{-2} r^2 
\int_{ \partial B(r)} \Vert (F_{A_1})_{\psi \psi} \Vert^2_g dV_g.
\end{align}
Combining the above, and letting $r_0 = r$, we have 
\begin{align*}
 \int_{D(r)} \Vert Rm_{g} \Vert^2_g dV_g
& \leq  c(r) \cdot r
 \left\{ \int_{ \partial B(r)} \Vert (R_{A_1})_{\psi \psi} \Vert^2_g dV_g
\right\}^{1/2}
\cdot 
 \left\{ \int_{\partial B(r)} \Vert (R_{A_1})_{r \psi} \Vert^2_g d \sigma_g 
\right\}^{1/2}\\
&+  C r^{- 2 \alpha}\\
& \leq \frac{ c(r)}{2} \cdot r
\int_{ \partial B(r)} \Vert Rm_g \Vert^2_g dV_g+  C r^{-2\alpha},
\end{align*}
where $c(r) =  (2 - \epsilon'(r))^{-1} ( 1 - C \epsilon(r) - \delta)^{-1}$.
Therefore for all $r$ sufficiently large, choosing $\delta$ 
sufficiently small, there exists a small constant $\delta' > 0$ 
\begin{align*}
 \int_{D(r)} \Vert Rm_{g} \Vert^2_g dV_g
& \leq  \frac{r}{4} ( 1 + \delta')
\int_{ \partial B(r)} \Vert Rm_g \Vert^2_g dV_g+  C r^{-2 \alpha}.
\end{align*}
We rewrite this inequality as 
\begin{align}
f(r) \leq c_1 r f'(r) + c_2 r^{- 2 \alpha},
\end{align}
where $c_1 =  - \frac{1}{4} ( 1 + \delta')$.
Since $\alpha < 2$, we have $0 < 1 + 2 c_1 \alpha$, so let 
\begin{align}
c_3 = \frac{c_2}{1 + 2 c_1 \alpha}, 
\end{align}
and consider the function $h(r) = \mbox{max} \{ f(r) - c_3 r^{- 2\alpha}, 0\}$.
When $h > 0$ we have 
\begin{align*}
c_1 r h'(r) &= c_1 r(f'(r) + 2 \alpha  c_3 r^{- 2\alpha-1})
\geq f(r) - c_2 r^{- 2\alpha } + 2c_1 c_3 \alpha r^{-2\alpha}\\
& = f(r) - c_3 r^{-2\alpha },
\end{align*}
and therefore $h(r)$ satisfies the inequality
\begin{align}
h(r) \leq c_1 r h'(r).
\end{align}
Integrating this inequality, we obtain $ h(r) \leq h(1) \cdot r^{1/c_1}$,
and therefore $f(r) \leq c_5 r^{ - 2\beta}$,
where $2\beta = \mbox{max} \{ -2 \alpha, 1/c_1 \}$.
The pointwise decay (\ref{pdc}) follows from (\ref{4}).
\end{proof}
To finish the proof, Lemma \ref{bali} and 
Theorem \ref{higherlocalregthm2} give improved 
pointwise decay on the full curvature tensor and its 
covariant derivatives, and the ALE property then follows 
by the result in \cite[Theorem 1.1]{BKN}.
\section{Constraints}
\label{Constraints}
For an ALE space $X$ with several ends, we have the 
signature formula
\begin{align}
\label{signature}
\tau(X) = \frac{1}{12 \pi^2} \left( \int_X |W^+_g|^2 dV_g - \int_X
|W^-_g|^2 dV_g \right) - \sum_i \eta( S^3 / \Gamma_i ),
\end{align}
where $\Gamma_i \subset SO(4)$ is the group corresponding to
the $i$th end, and $\eta(  S^3 / \Gamma_i)$ is the 
$\eta$-invariant. The Gauss-Bonnet formula in this context is
\begin{align}
\label{GB}
\chi(X) = \frac{1}{8 \pi^2} \left(
\int_X |W_g|^2 dV_g 
- \frac{1}{2} \int_X |Ric_g|^2dV_g + \frac{1}{6} \int_X R_g^2 dV_g \right)
+ \sum_i \frac{1}{|\Gamma_i|}.
\end{align}
See \cite{Hitchin} for a nice discussion of these formulas.
We remark that these put constraints on the ends of 
the ALE spaces that can arise in Theorem \ref{decayale}. 
For example, if $X$ is locally conformally flat and
scalar flat, then 
\begin{align}
\label{eta}
\tau(X) &= - \sum_i \eta( S^3 / \Gamma_i ),\\
\label{grb}
\chi(X) &= - \frac{1}{16 \pi^2} \int_X |Ric_g|^2dV_g 
+ \sum_i \frac{1}{|\Gamma_i|}.
\end{align}
Equation (\ref{eta}) says that the sum of the 
$\eta$-invariants must be integral, 
and (\ref{grb}) gives the inequality, 
\begin{align}
\chi(X) \leq  \sum_i \frac{1}{|\Gamma_i|} \leq \# \{ \mbox{ends of } X \}. 
\end{align}
If we assume that $X$ is scalar-flat and anti-self-dual we obtain 
\begin{align}
3 \tau(X) &= - \frac{1}{4 \pi^2} \int_X
|W^-_g|^2 dV_g  -3 \sum_i \eta( S^3 / \Gamma_i ),\\
2 \chi(X) &= \frac{1}{4 \pi^2} \left(
\int_X |W^-_g|^2 dV_g 
- \frac{1}{2} \int_X |Ric_g|^2dV_g \right) + 2 \sum_i \frac{1}{|\Gamma_i|}.
\end{align}
Adding these together,
\begin{align}
2 \chi(X) + 3 \tau(X) =   -  \frac{1}{8 \pi^2} 
\int_X |Ric_g|^2dV_g +  \sum_i 
\Big( \frac{2}{|\Gamma_i|} - 3 \eta( S^3 / \Gamma_i ) \Big),
\end{align}
which implies the inequality
\begin{align}
2 \chi(X) + 3 \tau(X)\leq  \sum_i 
\Big( \frac{2}{|\Gamma_i|} - 3 \eta( S^3 / \Gamma_i ) \Big),
\end{align}
with equality if and only if $g$ is Ricci-flat. 

We conclude this section with a pinching theorem.
\begin{theorem} 
\label{pinching}
Let $(M_1,g)$ be a complete, oriented, noncompact 4-dimensional Riemannian 
manifold with $g$ scalar-flat and anti-self-dual.
Assume that
\begin{align}
\ C_S < \infty, \mbox{ and }  b_1(M) = 0.
\end{align}
Then there exists an $\epsilon_1$ depending 
upon $C_S$ such that if $\Vert Rm \Vert_{L^2} < \epsilon_1$
then $(M_1,g)$ is isometric to $\mathbf{R}^4$. 
\end{theorem}
\begin{proof}
From  \cite[Proposition 1]{Carron1}), there exists
$\epsilon_1 > 0$ (depending on the Sobolev constant),
so that if  $\Vert Rm \Vert_{L^2} < \epsilon_1$,
then the first $L^2$ Betti number vanishes, 
which implies that $(M_1,g)$ has only one end. 
From Theorem  \ref{decayale}, this end is ALE of order $\tau <2$.
Poincare duality gives $b_3(M_1) = b_1(M_1)$, 
$b_4(M_1) = 0$, and the Gauss-Bonnet formula is 
\begin{align}
\label{GB3}
1 + b_2(M_1) - \frac{1}{| \tilde{\Gamma}|}
= \frac{1}{8 \pi^2} \int_{M_1} \left( - 
\frac{1}{2} |Ric|^2 + |W^-|^2 \right)dV_g.  
\end{align}
If the right hand side is sufficiently small, we conclude 
that $ \tilde{\Gamma} = \{ e \}$, and $b_2(M_1) = 0$.
Since $(M_1,g)$ has exactly one asymptotically 
Euclidean end and since $(M_1,g)$ is half-conformally flat, 
from the signature formula we conclude that $(M_1,g)$ is conformally flat. 
From (\ref{GB3}), $(M_1,g)$ is therefore also Ricci-flat. By Bishop's relative 
volume comparison theorem \cite{Bishop}, 
$(M_1,g)$ is isometric to Euclidean space.  
\end{proof}
\bibliography{BachFlatALE_references}
\vspace{-1mm}
\end{document}